\documentclass[11pt]{amsart}

\usepackage{graphicx,multirow,rotating}

\newtheorem{thm}{Theorem}[section]

\newtheorem{rem}{Remark}[section]

\def\rot#1{\begin{sideways}{#1}\end{sideways}}

\def\1{{{\mbox{${\rm{1\negthinspace\negthinspace I}}$}}}}
\newcommand{\eref}[1]{(\ref{#1})}

\newcommand\ind{{ {{1}}\hspace{-0,8mm}{\mathrm I}}}

\renewcommand{\baselinestretch}{1.25}
\newcounter{hypc}
\newcommand{\hypothese}[2]{\stepcounter{hypc}
\tag{$\mathbf{A}_{\thehypc}^{#2}$}
\label{#1}}
\newcounter{cond}
\newcommand{\condition}[1]{\stepcounter{cond}
\tag{$\mathbf{R}_{\thecond}^X$}
\label{#1}}

\begin{document}
\title[Finite sample penalization in adaptive density deconvolution]{Finite sample penalization
in adaptive density deconvolution.}
\author[F. Comte]{F. Comte$^{1}$}\thanks{$^1$ Universit\'e Paris V,  MAP5, UMR CNRS
8145.}
\author[Y. Rozenholc]{Y. Rozenholc$^{1}$}
\author[M.-L. Taupin]{M.-L. Taupin$^3$}\thanks{$^3$ IUT de Paris V et
Universit\'e d'Orsay, Laboratoire de Probabilit\'es, Statistique
et Mod\'elisation, UMR 8628.}
\begin{abstract}
We consider the problem of estimating the density $g$ of identically distributed variables $X_i$, from a sample $Z_1,
\dots, Z_n$ where $Z_i=X_i+\sigma\varepsilon_i$, $i=1, \dots, n$ and $\sigma
\varepsilon_i$ is a noise independent of $X_i$ with known density $
\sigma^{-1}f_\varepsilon(./\sigma)$. We generalize adaptive estimators,
constructed by a model selection procedure, described in Comte et al.~(2005).
We study numerically their properties in various contexts and we test their
robustness. Comparisons are made with respect to deconvolution kernel
estimators, misspecification of errors, dependency,... It appears that our
estimation algorithm, based on a fast procedure, performs very well in all
contexts.

\end{abstract}
\maketitle
\begin{center}
\today
\end{center}
\noindent {\sc {\bf Keywords.}} {\small  Adaptive estimation.
Density deconvolution. Model selection. Penalized contrast.
Projection estimator. Simulation study. Data-driven.} \\

\section{Introduction}
In this paper, we consider the problem of the nonparametric
density deconvolution of $g$, the density of identically distributed
variables $X_i$, from a sample $Z_1, \dots, Z_n$ in the model
\begin{eqnarray}
\label{model}
Z_i=X_i+\sigma \varepsilon_i, \;\; i=1, \dots, n,\end{eqnarray}
where the $X_i$'s and
$\varepsilon_i$'s are independent sequences, the $\varepsilon_i$'s
are i.i.d. centered random variables with common density $f_{\varepsilon}$,
that is $\sigma\varepsilon_i$ is a noise with known density $\sigma^{-1}
f_{\varepsilon}(./\sigma)$ and known noise level $\sigma$.

Due to the independence between the $X_i$'s and the $\varepsilon_i$'s,
the problem is to estimate $g$ using the observations $Z_1,\cdots,Z_n$ with
common density $f_Z(z)=\sigma^{-1}g\star f_\varepsilon(./\sigma)(z).$
The function $\sigma^{-1}f_{\varepsilon}(./\sigma)$ is often called the
convolution kernel and is completely known here.

Denoting by $u^*$ the Fourier transform of $u$, it is well known that since
$g^*(.)=f_Z^*(.)/f_\varepsilon^*(\sigma.)$,  two factors
determine the estimation accuracy in the standard density deconvolution problem : the
smoothness of the density to be estimated,
and the one of the error density which are described
by the rate of decay of their Fourier transforms.
In this context, two classes of errors are usually considered: first the so called ``ordinary
smooth" errors with polynomial decay of their  Fourier transform
and second, the ``super smooth" errors with Fourier transform having an exponential decay.

For further references about density deconvolution see e.g.
Carroll and Hall~(1988), Devroye (1989), Fan~(1991a, b), Liu and Taylor~(1989),
Masry~(1991, 1993a, b), Stefansky~(1990), Stefansky and
Carroll~(1990), Taylor and Zhang~(1990), Zhang~(1990) and Cator~(2001), Pensky and
Vidakovic~(1999), Pensky~(2002), Fan and Koo~(2002), Butucea~(2004), Butucea
and Tsybakov~(2004), Koo~(1999).

The aim of the present paper is to provide a complete simulation study of the
deconvolution estimator constructed by a penalized contrast minimization on a model $S_m$,
a space of square integrable functions having
a Fourier transform with compact support included into $[-\ell_m,\ell_m]$ with
$\ell_m=\pi L_m$.
Comte et al.~(2005) show that for $L_m$ being a positive integer, this penalized contrast
minimization selects the relevant projection space $S_m$ without any prior information
on the unknown
density $g$. In most cases, it is an adaptive estimator in the sense that it achieves the optimal rate of convergence
in the minimax sense, studied by Fan~(1991a), Butucea~(2004) and Butucea and Tsybakov~(2004).
It is noteworthy that, contrary to what usually happens, $\ell_m$ does not correspond here
to the dimension of the projection space but to the length of the support of the Fourier
transform of the functions of $S_{m}$. Thus we will refer in the following
to $\ell_m$ as the "length" of the model $S_{m}$.

Moreover, in the context of integer $L_m$, Comte et al.~(2005) provide a brief
simulation which shows that the selected $L_m$ are rather small
and therefore far from the asymptotic.
Our present study shows that it is relevant to choose $\ell_m=\pi L_m$ on a thinner
grid than one included in $\pi \mathbb{N}$.

Thus we start by stating a modification of the results in Comte et al.~(2005)
to take into account this thinner grid of values $\ell_m$
and we show that the resulting penalized minimum contrast estimator is
an adaptive estimator in the sense that it achieves the optimal rate of convergence
in the minimax sense.
Here, the penalty depends on the smoothness of the errors density and therefore
we consider two cases: Laplace density (ordinary smooth) and Gaussian density (super smooth).

We illustrate, through examples, the influence of over-penalization and
under-penalization and propose practical calibrations of the penalty in
all considered cases.

Then we study in very large simulations the  non asymptotic
properties of our estimator by considering various types of
densities $g$, with various smoothness properties like Cauchy distribution, Gaussian density
and finally F\'ejer-de-la-Vall\'ee Poussin-type density.

We present some examples, that illustrate how the algorithm works.
We give the mean integrated squared error (MISE) for the two types of errors density, for all
the test densities, for various $\sigma$, and for various sample size.
Our results present global tables of MISE and comparisons between MISE and the theoretical expected
rates of convergence.

Lastly, the robustness of our procedure is tested in various ways:
when the observations are dependent, when $\sigma$ is very small
(leading to a problem of density estimation) and when the errors
density $f_\varepsilon$ is misspecified or not taken into account.
In those cases, we compare our procedure with previous results of
Delaigle and Gjibels~(2004a, 2004b) and Dalelane~(2004) (direct
density estimation).

The conclusions of our study are the following.
Our estimation procedure provides very good results; better than the kernel
deconvolution methods described and studied in Delaigle
and Gijbels~(2004a). Our estimation procedure is robust when the $Z_i$'s
are no longer independent and even not strongly mixing.
We underline the importance of the noise level in the quality of
estimation, and we check that, in the case of a very small noise, we obtain
MISE's that have the same order as some recent results obtained by
Dalelane~(2004) for direct density estimation.
Lastly our results show that a misspecification of the errors density
slightly increases the error of estimation, but less than the use of the direct density estimator
(without deconvolving), as it was already mentioned in Hesse~(1999).
>From a practical point of view it is important to note that our algorithm is a
fast algorithm ($O(n\ln(n))$ operations) based on the Inverse Fast Fourier Transform (IFFT).

The paper is organized as follows. In section 2, we present the
model, the assumptions, the adaptive estimator and its expected rates of convergence.
In Section 3, we describe the implementation of the estimates (see \ref{descri})
and the computations of the
associated integrated squared errors (\ref{MISE}). Section 4 presents the chosen penalties
(see \ref{penalites}) and describes the
framework of our simulations.
The simulation results are gathered in Section 5 and an appendix is devoted to
the proof of our theorem.

\section{General framework and theoretical results}
\subsection{Notations and assumptions}\label{secmodel}

 For $u$ and $v$ two square integrable functions, we denote by
$u^*$ the Fourier transform of $u$, $u^*(x)=\int e^{itx}u(t)dt$
and by $u*v$ the convolution product, $u*v(x)=\int
u(y)v(x-y)dy$. Moreover, we denote by $\|u\|^2=\int_{{\mathbb R}} |u(x)|^2 dx$.

Consider Model \eref{model} under the following assumptions.
\begin{align}
& \mbox{ The } X_i\mbox{'s and the } \varepsilon_i\mbox{'s
are independent and identically distributed random} \hypothese{iid}{ }\\
&\mbox{ variables and the sequences }(X_i)_{i\in \mathbb{N}}\mbox{ and }
(\varepsilon_i)_{i\in \mathbb{N}} \mbox{ are independent.} \notag \\
&\mbox{ The density } f_\varepsilon \mbox{ belongs to
  }\mathbb{L}_2(\mathbb{R})  \mbox{ and is such that for all }
x\in \mathbb{R}, \hypothese{fepspair}{\varepsilon} f_\varepsilon^*(x)\not=0\notag.
\end{align}
Under assumption \eref{iid}, the $Z_i$'s are independent and identically distributed
random variables. Assumption \eref{fepspair}, usual for the
construction of an estimator in density deconvolution, ensures
that $g$ is identifiable.

The rate of convergence for estimating $g$ is strongly related to
the rate of decrease of the Fourier transform of the errors
density $f_\varepsilon^*(x)$  as $x$ goes to infinity. More
precisely, the smoother $f_\varepsilon$, the quicker the rate of decay of
$f^*_{\varepsilon}$ and the slower the rate of
convergence for estimating $g$. Indeed, if $f_\varepsilon$ is very
smooth, so is $f_Z$ the density of the observations $Z$ and thus it
is difficult to recover $g$.
This decrease of $f_\varepsilon^*$ is described by the following assumption.
\begin{align}
& \mbox{ There exist nonnegative real numbers } \gamma,~ \mu, \mbox{ and }\delta \mbox{ such
that }\hypothese{regufeps}{\varepsilon} \\
&\notag\vert f_\varepsilon^*(x)\vert\geq \kappa_0(x^2+1)^{-\gamma/ 2}\exp\{-\mu\vert x\vert^\delta\}
\end{align}
 When $\delta=0$ in assumption \eref{regufeps}, $f_\varepsilon$ is usually called ``ordinary smooth'', and
when  $\mu>0$ and $\delta>0$, the error density is usually called
``super smooth''. Indeed densities satisfying assumption \eref{regufeps} with
$\delta>0$ and $\mu>0$ are infinitely differentiable. For instance, Gaussian
or Cauchy distributions are super smooth of order $\gamma=0,
\delta=2$ and $\gamma=0, \delta=1$ respectively, and the symmetric exponential
(also called Laplace) distribution with $\delta=0=\mu$ and $\gamma=2$ is an
ordinary smooth density. Furthermore, when $\delta =0$,
\eref{fepspair} requires that $\gamma> 1/2$ in \eref{regufeps}.
By convention, we set $\mu=0$ when $\delta=0$ and we assume that
$\mu>0$ when $\delta> 0$. In the same way, if $\sigma=0$, the
$X_i$'s are directly observed without noise and we set
$\mu=\gamma=\delta =0$.

For the construction of the estimator we need the following more
technical assumption.
\begin{align}
&\mbox{The density } g\mbox{  belongs to }\mathbb{L}_2(\mathbb{R})\mbox{ and there
exists some positive real } M_2\hypothese{momg}{X}\\
&\mbox{ such that } g \mbox{  belongs to } \left\{t \mbox{ density such that }
  \int x^2t^2(x)dx \leq M_2 <\infty\right\}
\notag.\end{align}
This assumption \eref{momg}, quite unusual but unrestrictive,
  already appears in density deconvolution in a slightly different way
in  Pensky and Vidakovic~(1999) who assume, instead of \eref{momg}
that $\sup_{x \in \mathbb{R}} \vert x\vert g(x)<\infty$.  The main drawback of this
condition is that it is not stable by translation, but an empirical centering
of the data seems to avoid
practical problems. 

Since rates of convergence  depend on the smoothness of $g$ we
introduce regularity conditions.
\begin{align}\condition{super}
& \mbox{There exists some positive real numbers } s, r, b \mbox{ such that
   the density } \\ & g \in
\mathcal{ S}_{s,r,b}(C_1)=\left\{ t
\mbox{ density }\; : \; \int_{-\infty}^{+\infty}
|t^*(x)|^2(x^2+1)^{s}\exp\{2b |x|^{r}\} dx\leq
C_1 \right\}. \notag\\
& \mbox{There exists some positive real numbers } K \mbox{ and }d\mbox{
  such that the density } \condition{entiere} \\ &  g \in \mathcal{S}_d(C_2)=\left\{ t \mbox{ density such that
  for all }x \in \mathbb{R}, \;  \vert t^*(x)\vert \leq C_2\ind_{[-d,d]}(x)\right\}.\notag
\end{align}
Note that densities satisfying \eref{super} with
$r=0$ belong to some Sobolev class of order $s$, whereas densities satisfying
\eref{super} with
$r>0, b>0 $ are infinitely differentiable.
Moreover, such densities admit analytic continuation on a finite
width strip when $r=1$  and on the whole complex plane if
$r=2$.
The densities satisfying \eref{entiere}, often called entire functions, admit analytic
continuation on the whole complex plane (see Ibragimov and Hasminskii
~(1983)).

In order to clarify the notations, we denote by greek letters the parameters
related to the known distribution of the noise $\varepsilon$ and by latin
letters the parameters related to the unknown distribution $g$ of $X$.

Let us now present and motivate the estimator.
\newpage
\subsection{The projection spaces and the estimators}
\subsubsection{Projection spaces}
 Let
$\varphi(x)=\sin(\pi x)/(\pi x)$ and $\varphi_{m,j}(x) = \sqrt{{L_m}}
\varphi({L_m}x-j)$.
Using that $\{\varphi_{m,j}\}_{j \in \mathbb{Z}}$ is an orthonormal basis of
the space of square integrable functions having a Fourier
transform with compact support included into $[-\pi {L_m}, \pi {L_m}]=[-\ell_m, \ell_m]$ (see
Meyer~(1990)), we denote by
$S_{m}$ such a space and consider the collection of linear spaces
$(S_{m})_{m\in \mathcal M_n}$, with  $\ell_m=m\Delta$, $\Delta>0$, and $m \in \mathcal{M}_n$ with
${\mathcal M}_n=\{1, \dots, m_n\}$, as projection spaces. Consequently,
\begin{eqnarray*}
S_{m}= {\rm Vect}\{\varphi_{m,j}, \; j\in \mathbb{Z}\},\;
=\{f\in \mathbb{L}_2(\mathbb{R}), \mbox{ with } \mbox{supp}(f^*) \mbox{ included
  into }[-\ell_m, \ell_m ]\},
\end{eqnarray*}
and the orthogonal projection of $g$ on $S_{m}$, $g_{m}$
is given by $g_{m}=\sum_{j\in {\mathbb Z}} a_{m,j} \varphi_{m,j}$,
with $a_{m,j} = <\varphi_{m,j},g>$.
Since this orthogonal projection involves infinite sums, we consider in
practice, the truncated spaces $S_{m}^{(n)}$ defined as
$$S_{m}^{(n)}= {\rm Vect }\left\{\varphi_{m,j}, |j|\leq K_n\right\}$$ where
$K_n$ is an integer to be chosen later.
Associated to those spaces we consider the orthogonal projection of $g$ on $S_{m}^{(n)}$ denoted by
$g_{m}^{(n)}$ and given by
$g_{m}^{(n)}=\sum_{|j|\leq K_n} a_{m,j} \varphi_{m,j}
\mbox{ with } a_{m,j} = <\varphi_{m,j},g>.$
\subsubsection{The non penalized estimators}
Associate this collection of models to the following contrast
function, for $t$ belonging to some $S_{m}$ of the collection
$(S_{m})_{{L_m}\in {\mathcal M}_n}$
$$\gamma_n(t)= \|t\|^2
-\frac 2{n} \sum_{i=1}^nu_t^*(Z_i), \;\;\; \mbox{ with } \;\;\;u_t(x) = \frac 1{2
\pi} \left(\frac{t^*}{f_{\varepsilon}^*(\sigma .)}\right)(-x).$$ Since
$\mathbb{E}\left[u_t^*(Z_i)\right]
=\langle t, g\rangle,$ we
find that $\mathbb{E}(\gamma_n(t))=\|t-g\|^2 -\|g\|^2 $ which is minimum when
$t\equiv g$. Since $\gamma_n(t)$ estimates the $\mathbb{L}_2$ distance between
$t$ and $g$, it is well adapted for estimating $g$.
Associated to the collection of models, the collection of
the non penalized estimators $\hat g_{m}^{(n)}$ is defined by
\begin{equation}\label{tronque}  \hat g_{m}^{(n)} = \arg\min_{t\in S_{m}^{(n)}} \gamma_n(t).
\end{equation}
By using that $t\mapsto u_t$ is linear, and that $\{\varphi_{m,j}\}_{ \vert j\vert\leq K_n}$ is an
orthonormal basis of $S_{m}^{(n)}$, we have $\hat g_{m}^{(n)} =
\sum_{\vert j\vert \leq K_n} \hat a_{m,j} \varphi_{m,j} \;\;
\mbox{ where }\;\; \hat a_{m,j}= n^{-1}\sum_{i=1}^n
u_{\varphi_{m,j}}^*(Z_i),$ with $\mathbb{E}(\hat a_{m,j})=
<g,\varphi_{m,j}>=a_{m,j}.$

\subsubsection{The adaptive estimator}
The adaptive estimator is computed by using the
following penalized criteria
\begin{equation}\label{estitronc}
\tilde g=\hat g_{\hat m}^{(n)} \mbox{ with } \hat m= \arg\min_{m\in {\mathcal
    M}_n}
 \left[\gamma_n(\hat g_{m}^{(n)}) + \;
{\rm pen}(\ell_m)\right],
\end{equation}
where pen(.) is a penalty function based on the observations and the known distribution of
$\sigma\varepsilon_1$ without any prior information on $g$.

\subsection{Rate of convergence of the non adaptive estimator}\label{aim}
We recall here, using our setup, the bound for the risk of $\hat g_m$, proved in Comte et al.~(2005).
\begin{equation}\label{obj} {\mathbb E}(\|g-\hat g_{m}^{(n)}\|^2) \leq \|g-g_{m}\|^2
+ \| g_{m}-g_{m}^{(n)}\|^2 + \frac{{L_m}}{\pi n} \int
\left|\frac{\varphi^*(x)}{f^*_{\varepsilon}(\sigma{L_m}x)}\right|^2dx.\end{equation}

First, the variance term $$\frac{L_m}{\pi n}  \int
|\varphi^*(x)|^2|(f^*_{\varepsilon}(\sigma {L_m}x)|^{-2}dx=\frac{\ell_m}{\pi n}\int_{-1}^1
\frac{dx}{|f_{\varepsilon}^*(\sigma \ell_m x)|^2},$$  depends, as usual in
deconvolution problems, on the rate of decay of the Fourier
transform of $f_\varepsilon$, with larger variance for smoother
$f_\varepsilon$. Under assumption \eref{regufeps},  for $\ell_m\geq \ell_0$,
the variance term satisfies
$$\frac{\ell_m}{\pi n} \int_{-1}^1 \frac{dx}{\left|f^*_{\varepsilon}(\sigma
    \ell_m x)\right|^2} \leq \lambda_1\ell_m^{2\gamma+1-\delta}
\exp(2\mu(\sigma \ell_m)^{\delta})/n,$$
where
\begin{eqnarray}
\label{lambda1} \lambda_1= \frac{(\sigma^2+\ell_0^{-2})^{\gamma}}{\kappa_0^{2}R(\mu,\sigma,\delta)}\;  \mbox{ and }  R(\mu,\sigma, \delta)= \left\{
\begin{array}{ll} 1 & \mbox{ if } \delta =0 \\ 2\mu\delta\sigma^\delta & \mbox{
if } 0<\delta \leq 1 \\ 2\mu\sigma^\delta & \mbox{ if } \delta>1.\end{array}
\right.\end{eqnarray}

Second, under assumption \eref{momg}, $\|g_{m}-g_{m}^{(n)}\|^2$ is of order
$(M_2+1)\ell_m^2/(\pi^2K_n)$. Consequently, under
\eref{regufeps}, $K_n\geq (M_2+1)n$
ensures that the risk ${\mathbb E}(\|g-\hat g_{m}^{(n)}\|^2)$ has the  order
$$\parallel g-g_{m}\parallel^2+
 (2\lambda_1+1)\ell_m^{(2\gamma+1-\delta)}
\exp\left\{2\mu\sigma^\delta \ell_m^{\delta}\right\}/n.$$
Finally, the bias term $\|g-g_{m}\|^2$ depends on
the smoothness of the function $g$ and has the expected order for
classical smoothness classes since it is given by the distance
between $g$ and the classes of entire functions having Fourier
transform compactly supported on $[-\ell_m, \ell_m]$ (see
Ibragimov and Hasminskii~(1983)).


If $g$ satisfies \eref{entiere},
then the bias term $\parallel g-g_{m}\parallel^2=0$, by choosing $\ell_m=d$. It
follows that in that case the parametric rate of convergence for estimating $g$ is achieved.

If $g$ belongs to some ${\mathcal S}_{s,r,b}(C_1)$ defined by (\ref{super}), then the squared
bias term can be evaluated by using that
$$\|g-g_{m}\|^2=\frac 1{2\pi} \int_{|x|\geq \ell_m} |g^*(x)|^2 dx\leq
\frac{C_1}{2\pi}(\ell_m^2+1)^{-s}\exp\{-2b \ell_m^{r}\}.$$

Consequently, under \eref{momg}, if $K_n\geq (M_2+1)n$, the rate of convergence
of $\hat g_{m}^{(n)}$ is obtained by selecting the space $S_{m}^{(n)}$, and thus ${\ell_m}$, that
minimizes
$$\frac{C_1}{2\pi}(\ell_m^2+1)^{-s}\exp\{-2b \ell_m^{r}\} +
(2\lambda_1+1)\frac{\ell_m^{(2\gamma+1-\delta)}
\exp\left\{2\mu\sigma^\delta\ell_m^{\delta}\right\}}{n}.$$
One can see that if $\ell_m$ becomes too large, the risk explodes, due to the presence
of the second term.
Hence $\ell_m$ appears to be the  cut between the relevant low frequencies used
in the Fourier transforms to compute the estimate and the high frequencies which are
not used (and may even degrade the quality of the risk).

We give the resulting rates in Table \ref{rates}. For a density $g$ satisfying (\ref{super}),
rates are, in most cases, known to be the optimal one in the
minimax sense (see Fan~(1991a), Butucea~(2004), Butucea and
Tsybakov~(2004)). We refer to Comte et al.~(2005) for further
discussion about optimality.

\begin{center}
\begin{table}[!ht]
\begin{tabular}{|c|l||c|c||}\cline{3-4}\cline{3-4}
\multicolumn{2}{c||}{} &\multicolumn{2}{c||}{$f_\varepsilon$} \\\cline{3-4}
\multicolumn{2}{c||}{} & $\delta=0$ & $ \delta>0$ \\
\multicolumn{2}{c||}{} & ordinary smooth & super smooth \\\hline\hline
\multirow{4}{.2cm}{\\\vfill\null $g$} & $\;$ & $\;$ & $\;$ \\
& $\begin{array}{l}
  r=0\\
  \small{\mbox{Sobolev}(s)}
\end{array}$ &
$\begin{array}{l}
  \ell_{\breve m}=O(n^{1/(2s+2\gamma +1)})\\
  \mbox{rate}=O(n^{-2s/(2s+2\gamma+1)})\\
\mbox{{\it optimal rate}}
\end{array}$  &
$\begin{array}{l}
  \ell_{\breve m}=[\ln(n)/(2\mu\sigma^{\delta}+1)]^{1/\delta}\\
  \mbox{rate}=O( (\ln(n))^{-2s/\delta})\\
\mbox{{\it optimal rate}}
\end{array}$  \\
$\;$ & $\;$ & $\;$ & $\;$ \\
\cline{2-4}
& $\begin{array}{l}
  r>0\\
  \mathcal{C}^\infty
\end{array}$ &
$\begin{array}{l} \\
  \ell_{\breve m}=\left[{\ln(n)/2b}\right]^{1/r} \\
  \mbox{ rate}= \displaystyle  O\left(\frac{\ln(n)^{(2\gamma+1)/r}}n\right)\\
\mbox{{\it optimal rate}} \\
  \;\; \end{array}$ &
$\begin{array}{c}
  \ell_{\breve m}  \mbox{ implicit solution of } \\
  {\ell_{\breve m}}^{2s+2\gamma+1-r}\exp\{2\mu \sigma^\delta
 \ell_{\breve m}^\delta+2b \ell_{\breve m}^r\}\\
  \qquad= O(n)\\
\mbox{{\it optimal rate if }}r<\delta\\
\end{array}$
\\
\hline\hline
\end{tabular}\\~\\
 \caption{Optimal choice of the length ($\ell_{\breve m}$)
and resulting (optimal) rates.}\label{rates}
\end{table}
\end{center}
\vspace{-1cm}
In the case $\delta>0$, $r>0$, the rates are not explicitly given
in a general setting. For instance, if $r=\delta$, the rate is of order
\begin{equation}\label{vitextra} [\ln(n)]^b
n^{-b/(b+\mu\sigma^\delta)}  \mbox{ with } b=[-2s
\mu\sigma^\delta+(2\gamma -r+1)b]/[r (\mu\sigma^\delta+b)].
\end{equation} On the other hand, if $r/\delta\leq 1/2$, then
the rate is given by
\begin{equation}\label{vitforstable}
\ln(n)^{-2s/\delta} \exp\left[-2b
\left(\frac{\ln(n)}{2\mu\sigma^\delta}\right)^{r/\delta}\right].\end{equation}

\begin{rem}
\label{choixkn}
{\rm First, it is important to note that the condition $K_n\geq (M_2+1)n$ allows us to
 construct truncated spaces $S_{m}^{(n)}$ using $O(n)$ basis vectors and hence
to construct a tractable and fast algorithm from a practical point of view
(see Section 3).
Second, the choice of larger $K_n$ does not change the efficiency of our
estimator from a statistical point of view but only changes the speed of the
algorithm from a practical point of view.
}
\end{rem}

\subsection{Rate of convergence of the adaptive estimator}\label{secthm}

The following theorem is an extension of Theorems 4.1 and 4.2 in Comte et
al.~(2005).
This new version
states that, for any fixed $\Delta$, we can take $\ell_m=m\Delta$, with $m=1,\cdots, m_n$, instead of
$\ell_m=m\pi$.

\begin{thm}\label{genepenind}
Consider the model described in section \ref{secmodel} under
\eref{iid},\eref{fepspair},\eref{regufeps} and \eref{momg}
and the collection of estimators $\hat g_{m}^{(n)}$ defined by
(\ref{tronque}) with
$\ell_m=m\Delta$ for $m=1,\cdots, m_n$. Let
$\lambda_1$ and
$\lambda_2$ be two constants
depending on $\gamma,\kappa_0,\mu,\delta$ and $\sigma$.
Let $\kappa$ be some
numerical constant, not necessary the same in each case. Consider

\begin{tabular}{lll}
\textbf{1)}
&${\rm pen}(\ell_m)\geq \kappa \lambda_1
\ell_m^{2\gamma+1-\delta} \exp\{2\mu(\sigma \ell_m)^{\delta}\}/n$,
&if $0\leq \delta < 1/3$,\\

\textbf{2)}
&${\rm pen}(\ell_m)\geq \kappa [\lambda_1+
\mu\sigma^{1/3}\pi^{1/3}\lambda_2]
\ell_m^{2\gamma+2/3}\exp\{2\mu\sigma^{1/3}
\ell_m^{1/3}\}/n$,
&if $\delta=1/3$,\\

\textbf{3)}
&${\rm pen}(\ell_m)\geq\kappa [\lambda_1+ \mu\pi^{\delta}\lambda_2]
\ell_m^{2\gamma+((1/2+\delta/2)\wedge 1)}\exp\{2\mu (\sigma
\ell_m)^\delta\}/n$,
&if $\delta>1/3$,
\end{tabular}

\noindent then, if  $K_n\geq (M_2+1)n$ and $m_n$ is such that pen$(\ell_{m_n})$ is bounded,
the estimator $\tilde g=\hat
g_{\hat m}^{(n)}$ defined by \eref{estitronc} satisfies
\begin{equation}\label{resu1}
{\mathbb E}(\|g-\tilde g\|^2) \leq C\inf_{\ell_m\in \{1, \dots,
m_n\}}[\|g-g_{m}\|^2+ {\rm pen}(\ell_m)] + \frac{c}{\Delta n},
\end{equation} where $C$ and $c$ are constants depending on
$f_{\varepsilon}$.
\end{thm}

In the first two cases, the lower bound of the penalty
has the same order as the variance term and the risk of the adaptive
estimator $\tilde g$ has the order of the smallest risk among the estimators
associated to the collection of $\hat g_{m}^{(n)}$. Hence we get an adaptive to the smoothness of $g$
statistical procedure,  that can choose the optimal
$\ell_m$ in a purely data driven
way, up to the knowledge of $M_2$ through the choice of $K_n\geq (M_2+1)n$.

In the last case, a small loss  of order $\ell_m^{(3\delta/2-1/2)\wedge
\delta}$ may occur. Nevertheless, this loss does not affect the rate of
convergence if the bias is the dominating term, that is when $\delta >1/3$,
and $0<r<\delta$.
This loss changes the rate only when the variance is the dominating term, that is when
$1/3<\delta\leq r$ and consequently when the considered $\ell_m$ are powers of
$\ln(n))$. When $1/3<\delta\leq r$, the rate is faster than
logarithmic, and only a logarithm loss occurs, as a price to pay for adaptation.
This loss occurs in particular when both the density $g$ to be estimated and the
density of the errors $f_\varepsilon$ are
gaussian.

The interest of taking $\ell_m=m\Delta$ lies in the possibility of choosing the best
$\ell_m$ among more values. Nevertheless, the theorem highlights that too
small $\Delta$'s  make the remainder term $c/(n\Delta)$ become larger.
For instance, according to Table \ref{rates}, when $g$ satisfies \eref{super}, we can choose
$\Delta=1/\ln(n)$ and, when $\nu\leq 2$, since $\gamma>1/2$ (in order to guarantee that
$f_\varepsilon$ belongs to $\mathbb{L}_2(\mathbb{R})$), we do not lose anything in term of
rate of convergence. Clearly if $g$ is an entire function satisfying \eref{entiere},
$\Delta$ has to be fixed. Since we do not know in which smoothness class the true density
is, the only strategy ensuring that the good rate is achieved is to take a fixed $\Delta$.

\section{Estimates and associated MISE implementation}
\label{estimates}

\subsection{Steps of the simulations}
\noindent Given a density $g$, a distribution of error $\varepsilon$, a sample size $n$, a
value of $\sigma$, we sample the $Z_i$'s and do the following steps:

\noindent $-$ compute the estimators via their coefficients ($\hat a_{m,j}$).

\noindent $-$ compute the contrast using that
$$\gamma_n(\hat g_{m}^{(n)})=-\sum_{\vert j\vert\leq K_n} |\hat a_{m,j}|^2
 = -\|\hat g_{m}^{(n)}\|^2$$
\noindent $-$ minimize $\gamma_n(\hat
g_{m}^{(n)}) + {\rm pen}(\ell_m)$ and deduce the selected $\hat m$ and the
associated $\tilde g=\hat g^{(n)}_{\hat m}$

\noindent $-$ evaluate the estimation error by a computation of
the integrated squared error (ISE), $\|\tilde g-g\|^2$.

\noindent $-$ repeat all the previous steps 1000 times and compute an empirical
version of MISE, $\mathbb{E}\|\tilde g-g\|^2$.

\subsection{Computation of the estimators}\label{descri}
We fixed arbitrarily $\Delta=1/10$.
Given the data $Z_1, \dots, Z_n$, we need to compute for several
values of $\ell_m=\Delta,2\Delta, \dots$, the
coefficients of the estimate $\hat g_{m}^{(n)}$, $\hat
g_{m}^{(n)}=\sum_{\vert j\vert \leq K_n} \hat a_{m,j}
\varphi_{m,j}$, $\varphi_{m,j}=\sqrt{{L_m}}\varphi({L_m}x-j)$ with
$\varphi(x)=\sin(\pi x)/(\pi x)$. Since
$$\hat a_{m,j}= \frac 1n \sum_{k=1}^n
u^*_{\varphi_{m,j}}(Z_k) = \frac 1{2\pi n} \sum_{k=1}^n \int
e^{-ixZ_k} \frac{\varphi_{m,j}^*(x)}{f_{\varepsilon}^*(\sigma x)}dx$$ we get
that by denoting $\psi_Z(x) = n^{-1}\sum_{k=1}^n
e^{ixZ_k}$, the empirical Fourier transform of
$f_Z(.)=\sigma^{-1}g*f_{\varepsilon}(./\sigma))$,
then
\begin{eqnarray*}\hat a_{m,j} &=& \frac 1n \sum_{k=1}^n \frac
1{2\pi \sqrt{{L_m}}} \int_{-\pi {L_m}}^{\pi {L_m}}
\frac{e^{ix(Z_k-j/{L_m})}}{f_{\varepsilon}^*(\sigma x)} dx =
\frac{\sqrt{\ell_m}}{2\sqrt{\pi}} \int_{-1}^{1} e^{-2i\pi jx}\frac{\psi_Z(\ell_m x)
}{f_{\varepsilon}^*(\sigma \ell_m x)} dx.
\end{eqnarray*}
To compute integrals of type $2^{-1}
\int_{-1}^{1} e^{2i\pi jx} u(x)dx,$
we use their approximations via Riemann sums:
\begin{eqnarray}
\frac 1N \sum_{k=0}^{N-1} e^{ij\frac{-\pi+2k\pi}N}u(\frac{-1+ 2k}N).\label{FFT1}
\end{eqnarray}
Note that the IFFT (Inverse Fast Fourier Transform) Matlab function is
defined as the function which associates to a vector $(X(1),
\dots, X(N))'$ a vector $(Y(1), \dots, Y(N))'$ such that, for $N=2^M$,
\begin{eqnarray}
Y(j)= \frac 1N \sum_{k=1}^N X(k) e^{i(j-1)\frac{2\pi(k-1)}N} =
 \frac 1N \sum_{k=0}^{N-1} X(k+1) e^{i(j-1)\frac{2\pi k}N}.
\label{FFT2}
\end{eqnarray}

Hence, for
$X(k)=(\psi_Z/f_{\varepsilon}^*(\sigma .))(2(k-1) \ell_m/N)$ for $k=1, \dots, N$ and for
$Y=(Y_1, \dots, Y_N)'=IFFT(X)$, we get $\hat a_{m,j}= Y_{j+1}\sqrt{\ell_m/\pi}$ for
$j=0, \dots, N-1=2^M-1$. The quantity to be chosen is $M$ such that $K_n
=2^M-1\geq (M_2+1)n$. Indeed the $\hat a_{m,j}$'s can be
computed by using this IFFT with $K_n=N=2^M-1$ and with adequate shifts. In that way,
he quantity $\|g_{m}-g_{m}^{(n)}\|^2$
is always negligible with respect to the others.

One should take $M\geq \log_2(n+1)$.
After checking that a choice of a larger values (up to 11) does not change the estimation
quality, we finally choose $M=8$.

\subsection{Computation of the integrated squared error (ISE), $\|\tilde g-g\|^2$}
\label{MISE}
 We have two different ways for computing the integrated squared error
 $\|\tilde g-g\|^2$.
\begin{enumerate}
\item[(E1)] Standard approximation and discretization of
  the integral on an interval of ${\mathbb R}$ as it is done in Delaigle and
  Gijbels~(2004a) and Dalelane~(2004). In order to compare our results to
  theirs, we proceed to this valuation
on the same intervals.
\end{enumerate}
Since this evaluation on finite interval may lead to an
under-valuation of the ISE, we also propose an exact calculation of the
ISE on $\mathbb{R}$ as described in the following.
\begin{enumerate}

\item[(E2)]  Evaluation of the ISE on the whole real line.
We use the decomposition
\begin{eqnarray*}\|\hat g_{m}^{(n)}- g\|^2
&=&
\|g-g_{m}\|^2 +\|g_{m}-g_{m}^{(n)}\|^2 +\|g_{m}^{(n)}-\hat g_{m}^{(n)}\|^2.\end{eqnarray*}
In the cases we consider, $g^*$ is available and the bias term is computed by using the
standard formula $\|g-g_{m}\|^2 =
(1/(2\pi) \int_{|x|\geq \ell_m} |g^*(x)|^2dx$. We bound $\|g_{m}-g_{m}^{(n)}\|^2$ by a term of
order $\ell_m^2/K_n\leq \ell_m^2/2^M$.
Finally, the variance term $\|g_{m}-\hat g_{m}^{(n)}\|^2$, is calculated
using that
$$\|g_{m}^{(n)}-\hat g_{m}\|^2=\sum_{|j|\leq K_n} |a_{m,j}-\hat a_{m,j}|^2.$$
Consequently, we need the computation of $a_{m,j}= \sqrt{\ell_m}/(2\sqrt{\pi})\int_{-1}^{1} e^{-2\pi ijx}
g^*(\ell_m x)dx$, coefficients of the development of the projection
$g_{m}^{(n)}=\sum_{|j|\leq K_n}
a_{m,j}\varphi_{m,j}$ on $S_{m}^{(n)}$. Again, using  IFFT (see \eref{FFT1} and \eref{FFT2}), with $G=(G_1,\cdots,G_N)$ and
$G_k=g^*(2(k-1)\pi/N)$ for $k=1, \dots, N$, we get
$G^\star=(G_1^\star, \dots, G_N^\star)'=IFFT(G)$. Then $a_{m,j}=\sqrt{\ell_m/\pi}G^\star_{j+1}$
for $j=0, \dots, N-1=2^M-1$.
This second method requires the knowledge of $g^*$ and is unavoidable for
stable distributions for which the analytical form of $g$ is not available.
\end{enumerate}

\begin{rem}{\rm
\textbf{Speed of the algorithm:} Since the IFFT is a
fast algorithm, the computation of our estimates is also a fast algorithm and requires
only $O(2^M\ln(2^M))=O(n\ln(n))$ operations if $K_n=2^M-1$ is of order $n$. }

\end{rem}

\section{The practical framework}
\subsection{Description of the test densities $g$}\label{lesdensites}

We consider several types of densities $g$, and for each density, we give the interval $I$
on which the ISE is computed by the method (E1), which is the case
in all examples except for stable distributions, where the use of method (E2) is unavoidable.
The set of test densities can be split in three subsets.
First we consider densities having classical smoothness properties like H\"olderian
smoothness with polynomial decay of their Fourier transform. Second we consider
densities having stronger smoothness properties, with exponential decay of the
Fourier transform. And finally we consider densities with Fourier transform
compactly supported, that is satisfying Condition \eref{entiere}.

Except in the case of densities leading to infinite variance, we consider density
functions $g$ normalized with unit variance so that $1/\sigma^2$
represents the usual signal-to-noise ratio (variance of the signal
divided by the variance of the noise) and is denoted in the sequel
by $s2n$ defined as $s2n = 1/{\sigma^2}.$

\begin{itemize}
\item[{\bf (a)}] {\it Uniform distribution}: $g(x)= 1/(2\sqrt{3}) \1_{[-\sqrt{3},\sqrt{3}]}(x)$,
$g^*(x)=\sin(x\sqrt{3})/(x \sqrt{3}),$  $I=[-5,5]$.

\item[\textbf{(b)}] {\it Exponential distribution}: $g(x)=e^{-x}\1_{{\mathbb R}^+}(x)$, $g^*(x)= 1/(1-ix)$, $I=[-5,10]$.

\item[\textbf{(c)}] {\it $\chi^2(3)$-type distribution}: $X=1/\sqrt{6}U$, $g_X(x)=\sqrt{6}g(\sqrt{6}x)$, $U\sim \chi^2(3)$
where we know that $U\sim \Gamma(\frac 32, \frac 12)$,
$$g_U(x) = \frac 1{2^{5/2}\Gamma(3/2)} e^{-|x|/2}|x|^{1/2}, g_U^*(x)=
\frac 1{(1-2ix)^{3/2}},$$ and $I=[-1,16]$.

\item[\textbf{(d)}] {\it Laplace distribution}: as given in (\ref{doublexp}), $I=[-5,5]$.

\item[\textbf{(e)}] {\it Gamma distribution}: $\Gamma(2,3/2)$, with density
$g(x)=(3/2)^2x\exp(-3x/2)\1_{{\mathbb R}^+}(x)$,
$g^*(x)=-9/(4x^2+12ix-9)$. This density has variance 8/9, and is
renormalized for simulation, $I=[-5,25]$.

\item[\textbf{(f)}] {\it Mixed Gamma distribution}: $X=1/\sqrt{5.48} W$ with $W\sim 0.4 \Gamma(5,1) + 0.6 \Gamma(13,1)$,
$$g_W(x)=[0.4*\frac{x^4 e^{-x}}{\Gamma(5)} +
0.6\frac{x^{12}e^{-x}}{\Gamma(13)}]\mbox{1}\!\!\mbox{I}_{{\mathbb
R}^+}(x), g_W^*(x) = \frac{0.4}{(1-ix)^5} +
\frac{0.6}{(1-ix)^{13}},$$ and $I=[-1.5,26]$.

\item[\textbf{(g, h, i)}] {\it Stable distributions} of index $r=1/4$ (g), $r= 1/2$ (h),
$r=3/4$ (i). In those cases, the explicit form of $g$ is not available but we use that
$|g^*(x)|=\exp(-|x|^{r})$. The ISE is computed with method (E2).

\item[\textbf{(j)}] {\it Cauchy distribution}: $g(x)=(1/\pi)(1/(1+x^2))$, $g^*(x)=e^{-|x|}$, $I=[-10,10]$.

\item[\textbf{(k)}] {\it Gaussian distribution}: $X\sim {\mathcal N}(0,\sigma^2)$ with $\sigma=1$, $I=[-4,4]$.

\item[\textbf{(l)}] {\it Mixed Gaussian distribution}: $X\sim \sqrt{2}V$ with $V\sim 0.5{\mathcal N}(-3,1) + 0.5 {\mathcal N}(2,1)$
$$g_V(x) = 0.5 \frac 1{\sqrt{2\pi}}( e^{-(x+3)^2/2} + e^{-(x-2)^2/2)}),
\; g_V^*(x) = 0.5(e^{-3ix} + e^{2ix})e^{-x^2/2},$$ and $I=[-8,7]$.

\item[\textbf{(m,n, o, p)}] Scale transforms of the {\it F\'ejer-de la Vall\'ee-Poussin distribution}:
$$g(x)= \frac {1-\cos(px)}{p\pi x^2}, \;\; g^*(x)=(1-|x|/p)_+,$$ for $p=1$ in
\textbf{(m)}, $p=5$ in \textbf{(n)},  $p=10$ in \textbf{(o)} and $p=13$ in \textbf{(p)} and $I=[-10,10]$.
\end{itemize}

Densities \textbf{(a,b,c,d,e,f)} correspond to cases with $r=0$ (Sobolev
smoothness properties) with
different values of $s$, whereas densities \textbf{(g,h,i,j,k,l)} correspond to cases with $r>0$
(infinitely times differentiable) with different values for the power $r$.
Clearly, \textbf{(a,b)} are not even continuous.

Since the stable distributions \textbf{(g,h,i)} as well as the Cauchy
distribution \textbf{(j)}, have infinite variance, $s2n=1/\sigma^2$ is not
properly defined.

The stable distributions \textbf{(g,h,i)} also allow to study the
robustness of the estimation procedure when assumption \eref{momg}
is not fulfilled. When the density to be estimated $g$ is of type
\textbf{(g,h,i)} the tails of $g(x)$ are known to behave like
$|x|^{-(r+1)}$ (see Devroye (1986)).  It follows that, for such
densities, assumption \eref{momg} is fulfilled only if $r>1/2$.
Consequently only the stable distribution \textbf{(i)}, satisfies
\eref{momg}

The case of distributions \textbf{(m,n,o,p)} deserves some special comments: they correspond
to densities whose Fourier transform
has compact support included in $[-1,1]$ for \textbf{(m)}, $[-5,5]$ for \textbf{(n)}, $[-10,10]$
for \textbf{(o)} and $[-13,13]$ for \textbf{(p)}. As a consequence, the bias term $\int_{|x|\geq
\ell_m}|g^*(x)|^2dx$ equals zero as soon as $\ell_m\geq 1$
for \textbf{(m)}, $\ell_m\geq 5$ for \textbf{(n)}, for $\ell_m\geq 10$ for \textbf{(o)}, $\ell_m\geq 13$
for \textbf{(p)}. Therefore, the asymptotic rate for estimating this type of density is the parametric rate.

\begin{figure}[!ht]
$$\includegraphics[width=1\textwidth,
height=18cm,]{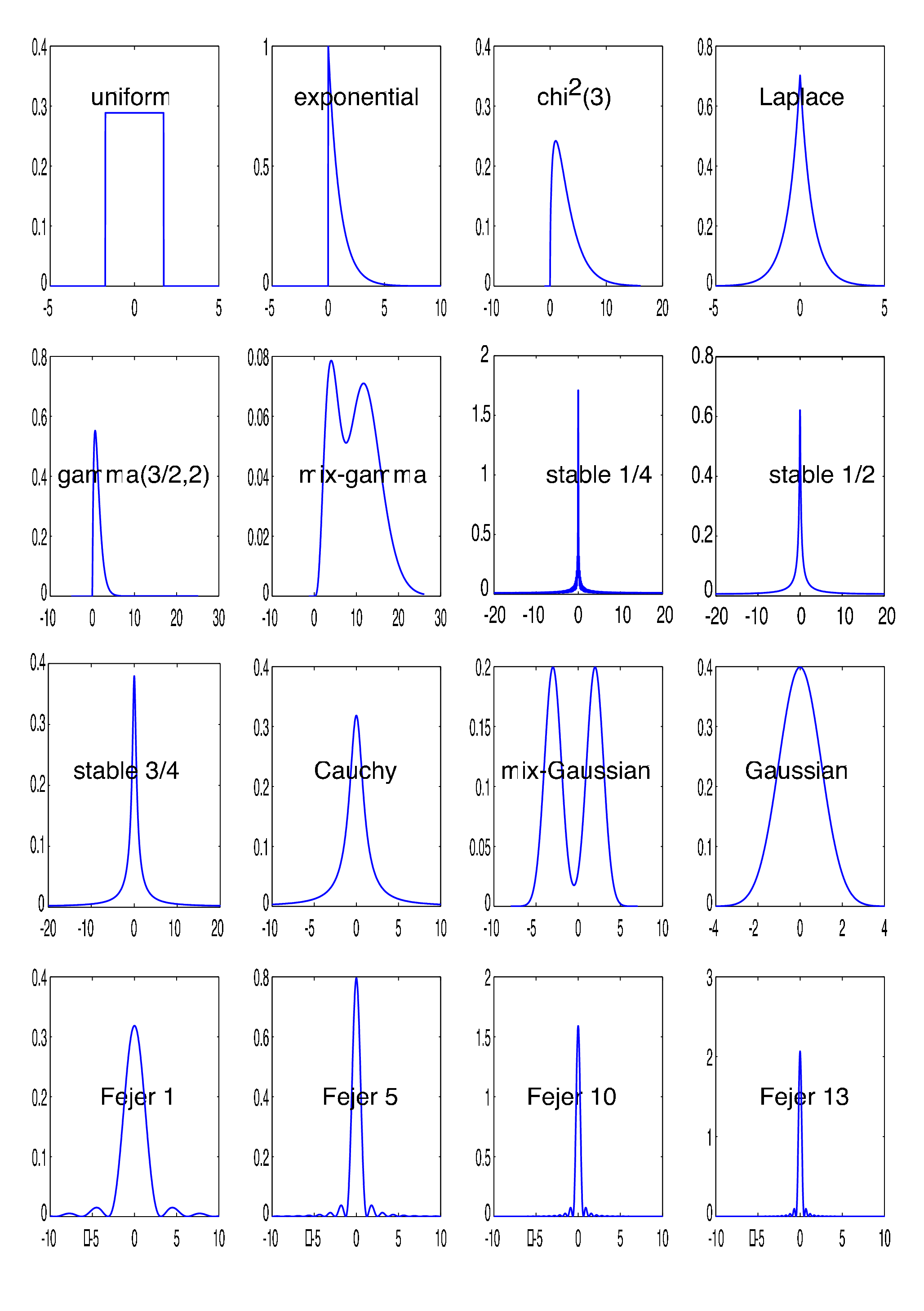}$$ \vspace{-1.3cm}\caption{Test
densities.}\label{functions}
\end{figure}

All above listed densities are plotted in Figure \ref{functions}. Note that for the
stable distributions, since no explicit form is available, we give in fact the
plot of the projection of the distribution on the space $S^{(n)}_{m}$ (for $\ell_m=10\pi$) as
computed by the projection algorithm.

We refer
to Devroye~(1986) for simulation algorithms of stable and Fejer-de la Vall\'ee-Poussin distributions.

\subsection{Two settings for the errors and the associated penalties.}
\label{penalites}

We consider two types of error density $f_\varepsilon$, the first one is the
Laplace distribution which is
ordinary smooth ($\delta=0$ in \eref{regufeps}), and the second one is the
Gaussian distribution which is super smooth
($\delta>0$ in \eref{regufeps}).

The penalty is connected to the variance order. In both settings, we will
precise this variance order and the value of the integral appearing in it. Since the
theory only gives the order of the penalty, by simulation
experiments, we fixed the constant
$\kappa$ and precise some additional negligible (with respect to the theory) terms
used to improve the practical results. In both cases we give the penalty given
in Comte et al. ~(2005) with $\Delta=\pi$ in $\ell_m=\Delta m$ and the new penalty allowing to
use a thinner grid for the $\ell_m$'s: here we take $\Delta=1/10$.
~\\

$\bullet$ \textbf{Case 1: Double exponential (or Laplace) $\varepsilon$'s.}

In this case, the  density of $\varepsilon$ is given by
\begin{equation}\label{doublexp}
f_{\varepsilon}(x)= e^{-\sqrt{2}|x|}/{\sqrt{2}}, \;\ f^*_{\varepsilon}(x)= (1+
x^2/2)^{-1}.\end{equation} This density corresponds to centered $\varepsilon$'s
with variance 1, and satisfying \eref{regufeps} with $\gamma=2$,
$\kappa_0=1/2$ and $\mu=\delta=0$.

The variance order is evaluated as
$$\kappa (\ell_m/(2\pi n))\int_{-1}^{1} 1/|f_{\varepsilon}^*(\sigma
\ell_m x)|^2dx=
\kappa (\ell_m/(\pi n)) \left(1+ \frac{\sigma^2 \ell_m^2}3 + \frac{\sigma^4\ell_m^4}{20} \right).$$
Let us recall that, in Comte et al.~(2005), $\Delta=\pi$, $\kappa=6\pi$ and the penalty is the following
\begin{eqnarray}
\label{oldpenlap}
{\rm pen}(\ell_m) = \frac{6}n\left[\ell_m+ \pi \ln^{2.5}(\ell_m/\pi) + \frac{\sigma^2\ell_m^3}3
 + \frac{\sigma^4 \ell_m^5}{20} \right].
\end{eqnarray}
The additional term $(\ln(\ell_m/\pi))^{2.5}$ is motivated by the works of Birg\'e
and Rozenholc~(2002) and Comte and Rozenholc~(2004). This term improves the quality of
the results by making the
penalty slightly heavier when $\ell_m$ becomes smaller.

Here, using intensive simulations study we propose the following penalty:
\begin{equation}\label{penlap} \fbox{\mbox{$\displaystyle{{\rm pen}(\ell_m) = \frac{2.5}n \left(1-\frac
1{s2n}\right)^2 \left[ \ell_m+ 8\ln^{2.5}(\zeta(\ell_m)) +
2\frac{\sigma^2\ell_m^3}3 + 3\left(1+ \frac 1{s2n}\right)^2
\frac{\sigma^4\ell_m^5}{10}\right],}$}}\end{equation}
with
\begin{equation}\label{zeta} \zeta(\ell_m)=\pi \1_{\ell_m<4}+\frac
{(\ell_m -2)^2}{4(\pi-2)}\1_{2\leq \ell_m<4}+ \ell_m \1_{\ell_m\geq
4}.\end{equation}

$\bullet$ \textbf{Case 2: Gaussian $\varepsilon$'s.} In that case, the errors density
$f_\varepsilon$ is given by
\begin{equation}\label{gaussian}
f_{\varepsilon}(x)=\frac 1{\sqrt{2\pi}}e^{-x^2/2}, \;\ f_{\varepsilon}^*(x)=
e^{-x^2/2}.\end{equation} This density satisfies \eref{regufeps} with $\gamma=0$,
$\kappa_0=1$, $\delta=2$ and $\mu=1/2$.

According to Theorem \ref{genepenind}, the penalty is slightly
heavier than the variance term, that is of order
$$\kappa \ell_m^{(3\delta/2-1/2)\wedge \delta}(\ell_m/(2\pi n))
\int_{-1}^{1} 1/|f_{\varepsilon}^*(\sigma \ell_m x)|^2dx=\kappa \ell_m^{(3\delta/2-1/2)\wedge
\delta}(\ell_m/(2\pi n))
\int_{-1}^{1} \exp(\sigma^2 \ell_m^2 x^2) dx.$$
Comte et al.~(2005), for $\Delta=\pi$,  choose $\kappa=6\pi$ and their penalty is the following
\begin{eqnarray}
\label{oldpengaus}
{\rm pen}(\ell_m)= \frac{6}n \left[\ell_m+\pi \ln^{2.5}(\ell_m/\pi) +
\frac{\ell_m^3\sigma^2}3\right]\int_{0}^{1} \exp[(\sigma\ell_m x)^2] dx.
\end{eqnarray}
According to the
theory, the loss, due to the adaptation is the term
$\sigma^2 \ell_m^2/3$. As previously, the additional term
$\ln(\ell_m/\pi)^{2.5}$ is motivated by simulations and the
works of  Birg\'e and Rozenholc~(2002) and Comte and
Rozenholc~(2004).

Using intensive simulation study we propose the following penalty
\begin{eqnarray}
\label{pengaus}
\fbox{\mbox{$\displaystyle{
{\rm pen}(\ell_m)= \frac{2.5}n \left(1-\frac 1{s2n}\right)^2 \left[\ell_m +
8\ln^{2.5}(\zeta(\ell_m)) +
\frac{\sigma^2\ell_m^3}3\right]\int_{0}^{1} \exp[(\sigma \ell_m x)^2] dx,}$}}\end{eqnarray}
where $\zeta(\ell_m)$ is defined by (\ref{zeta}) and the integral is
numerically computed.

\begin{figure}
$$\includegraphics[width=0.5\textwidth,height=7cm,]{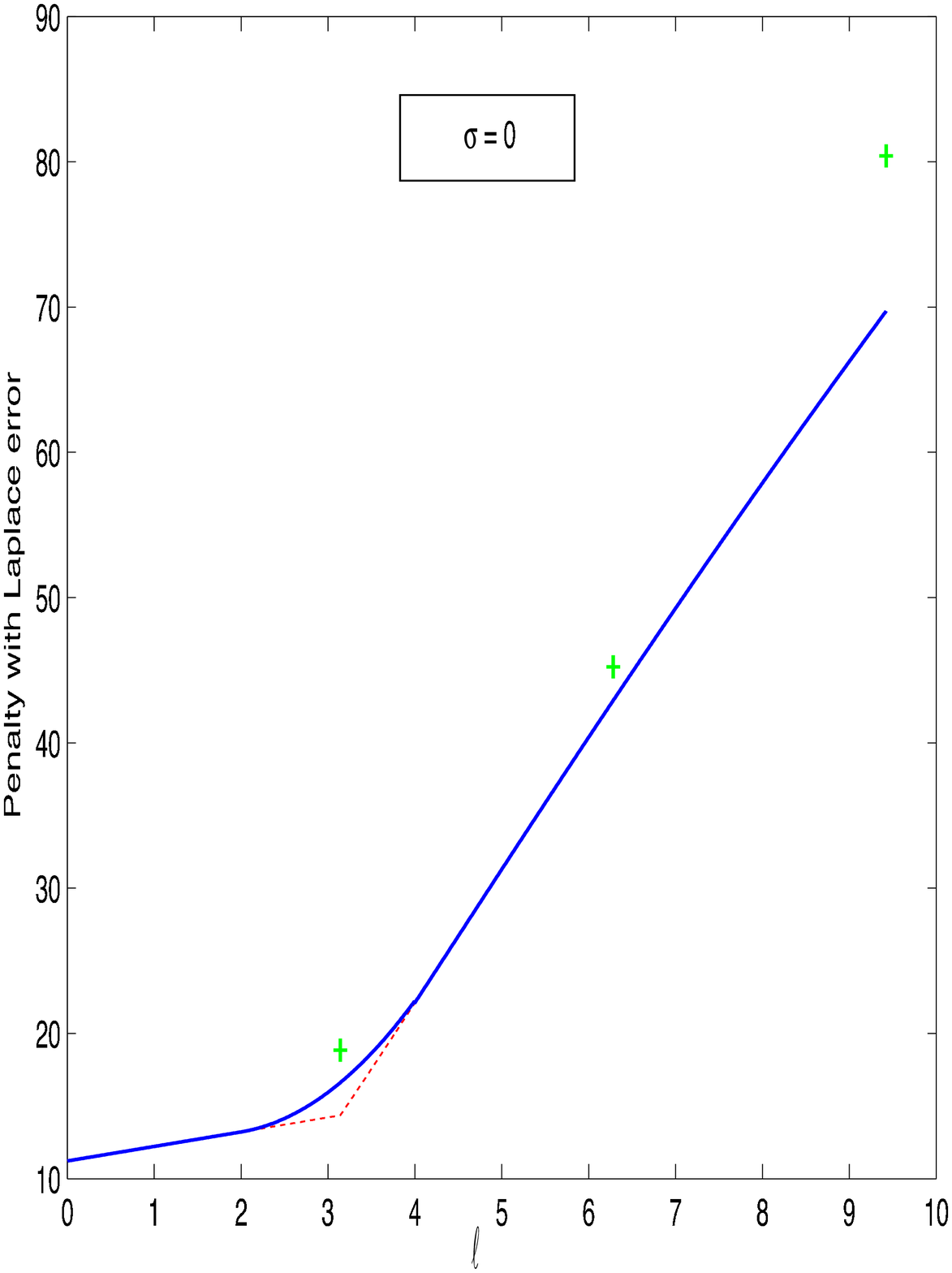}\includegraphics[width=0.5\textwidth,height=7cm,]{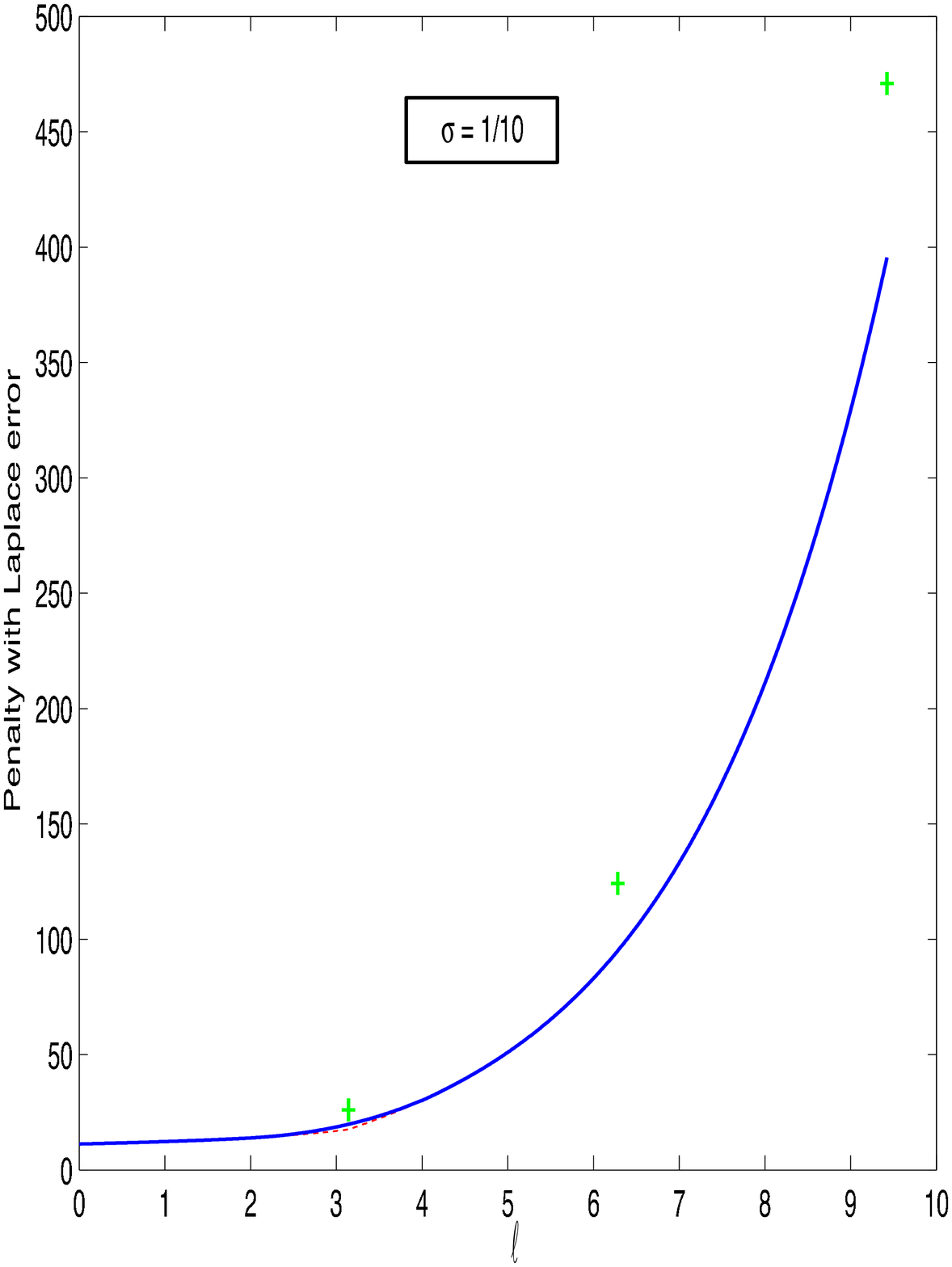}$$
\vspace{-1.3cm}\caption{Penalty in Comte et al.~2005 (crosses), penalty given by (\ref{penlap}) with $\zeta(\ell_m)=\pi\vee \ell_m$
(dotted line) and with $\zeta$ given by
(\ref{zeta}) (full line),  $\sigma^2=0$ and 0.1.}\label{penLs2=0}
\end{figure}

\begin{rem}
{\rm Note that when $\sigma=0$, both penalties are equal to $(2.5/n)(\ell_m+
8\ln(\zeta(\ell_m))^{2.5})$.}\end{rem}
\begin{rem}
{\rm
Since $\Delta=1/10$ we choose new constants and add a factor depending on $s2n$ in
\eref{penlap} and \eref{pengaus} with respect
to \eref{oldpenlap} and \eref{oldpengaus}. The
function $\zeta(\ell_m)$ is only chosen to give a smoother version
of $\ell_m\vee \pi$. The comparison of the penalty \eref{oldpenlap} for integer $L_m$'s,
the new penalty with $\zeta(\ell_m)=\ell_m\vee \pi$ (not smoothed) and
our final choice in \eref{penlap} is given in Figure \ref{penLs2=0} for $\sigma^2=0$
and for $\sigma^2=0.1$. The difference between the
two $\zeta$ functions clearly vanishes when $\sigma^2$ increases. }
\end{rem}

\begin{figure}
$$\includegraphics[width=0.7\textwidth,height=14cm, angle=-90]{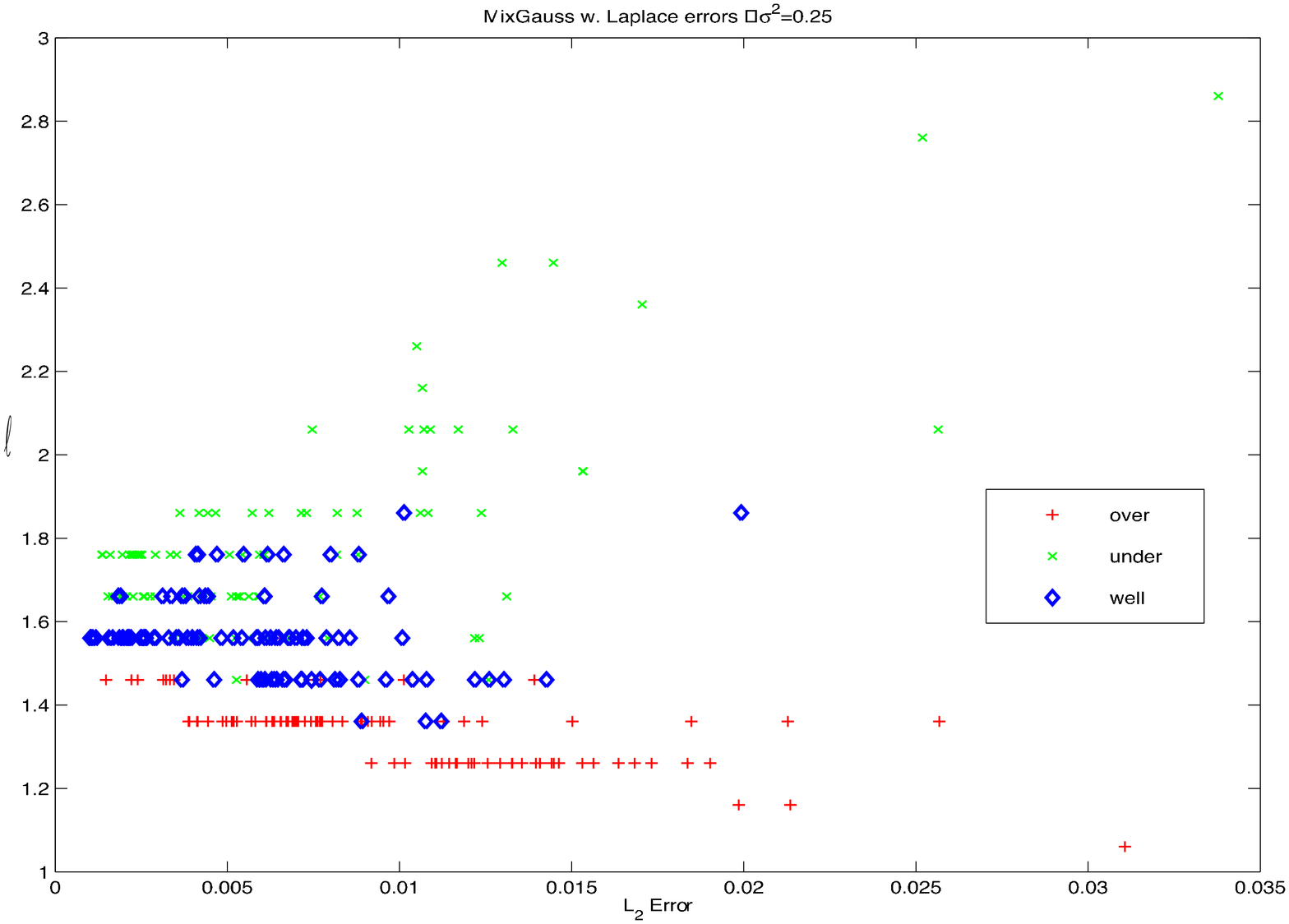}$$
\vspace{-0.8cm}\caption{Selected lengthes $\ell_m$ in function of the ISE ($L_2$
error) for the estimation of a mixed gaussian density in case of
under, over and good penalization.}\label{lmclouds}
\end{figure}

\begin{rem}{\rm

The influence of over- or under-penalization is illustrated in
Figure \ref{lmclouds}, where three penalties are tested for the
estimation of the mixed gaussian distribution. The figure plots the
selected $\ell_m$'s related to the ISE for 100 simulated path of
the distribution. This shows that over-penalization leads to smaller
selected $\ell_m$'s with increased ISE's, whereas
under-penalization leads to greater selected $\ell_m$'s with a more
important increase of both the dimensions and the ISE's.
The central cloud
of diamonds gives the selected $\ell_m$'s for our penalization and
shows that for this distribution our penalty is very well calibrated.

As illustrated by Figure \ref{lmclouds}, usually under penalization leads to larger values of
$\ell_m$ and increases the variance which degrades the MISE
more than over penalization. Hence it is better
to prevent from under penalization, the penalty is therefore increased.
Here, since $\ell_m$ takes
values on a thin grid, preventing against under penalization is less important
 and one can choose a smaller penalty which leads to a better trade-off
 between bias and variance. This leads to a better control of the risk. }
\end{rem}

\begin{rem}
{\rm
It is noteworthy that the penalty functions
\eref{penlap} and \eref{pengaus}  depend on
$s2n$ which is unknown. In Section \ref{robust}, we propose a study of the robustness
of the algorithm when
$s2n={\rm Var}(X)/\sigma^2={\rm Var}(Z)/\sigma^2-1$ is replaced by
a simple estimator (empirical variance of the observed $Z_i$'s instead
of the theoretical one). }\end{rem}

\subsection{Theoretical rates in our examples}

In order to compare the MISE resulting from our simulations, we
give in the Table \ref{theorates} the expected theoretical (and
asymptotic) rates corresponding to each cases we study.

It is noteworthy that even if theoretical results are established for
densities satisfying Condition \eref{super}, since we are in a simulation study, we
consider the explicit form of the Fourier transform of $g$ to evaluate the bias.
Consequently, for the calculation of the expected
theoretical rates given in Table \ref{theorates}, we denote by $s, r$ and $b$, the constants such that
\begin{equation}\label{convention}
\parallel g-g_{m}\parallel^2\leq \frac{1}{2\pi}\int_{\vert x\vert \geq \ell_m}\vert
g^*(x)\vert^2dx\leq \frac{A_s}{2\pi}(\ell_m^2+1)^{-s}\exp\{-2b
\ell_m^r\}.\end{equation}
Then, we evaluate the
theoretical rate of convergence by using the results in Table \ref{rates} with those $s, r$
and $b$.

Let us briefly comment this table \ref{theorates}.
Let us mention that with those choices of test densities, we describe all types of behavior of the
rates. According to Theorem \ref{genepenind}, except in the case where
$f_\varepsilon$ is the Gaussian density and the density to be estimated is
also the Gaussian density ($0\leq\delta\leq
1/3$ or $r<\delta$), the expected rates of convergence of the adaptive
estimator $\tilde g$ is the expected rate of convergence of the non penalized
estimator $\hat g_{\breve m}$ with  asymptotically optimal rate, that is the rate given
in Table \ref{rates}, with the convention \eref{convention}
about $s$, $r$ and $b$.

In the remainder case, when $f_\varepsilon$ is the Gaussian density and the density $g$ is
also the Gaussian density, $r=\delta=2>1/3$, the penalty is larger,  of a
logarithmic factor,
than the variance of the non penalized estimator $\hat g_{\breve m}$. Since
the penalty is the dominating term
in the trade-off with the bias, the rate of convergence of
$\tilde g$ is slower than the rate of convergence of the corresponding
non penalized estimator $\hat g_{\breve m}$. Let us be more precise. When $g$
is Gaussian, we have a bias term given by
$$\int_{|x|\geq \ell_m}|g^*(x)|^2dx = 2 \int_{\ell_m}^{+\infty} \exp(-x^2)dx\leq
2\int_{\ell_m}^{+\infty} \exp(-\ell_m x)dx\leq \frac{\exp(-\ell_m^2)}{\ell_m}$$
and a variance term of order ${\ell_m}^{-1}\exp(2\mu(\sigma\ell_m)^2)$.
So that, according to the convention \eref{convention}, we apply Formula (\ref{vitextra})
with $s=1/2$, $b=1/2$, $r=2$, $\delta=2$ and $\mu=1/2$, to get that the rate of convergence of the
non penalized estimator $\hat g_{\breve m}$ is of order
$$\ln(n)^{-\frac{1}{2}} n^{-\frac{1}{\sigma^2+1}}.$$
Now, according to Theorem \ref{genepenind}, the penalty is of order $\ell_m\exp(2\mu
(\sigma\ell_m)^2)$. We obtain that the rate of
convergence of the adaptive estimator $\tilde{g}$ is of order
$$(\ln(n))^{-\frac12\frac{\sigma^2 -1}{\sigma^2+1}}n^{-\frac 1{\sigma^2+1}}.$$
This implies a negligible loss of order $\ln(n)^{1/(1+\sigma^2)}$ for not
knowing the smoothness of $g$.

\begin{rem}
{\rm
Let us mention that taking $\sigma=0$ in columns 2 and 3 in Table \ref{theorates} does not always
provide the theoretical rates in the last column, with $\sigma=0$.
Some of the results above are not continuous when
$\sigma\rightarrow 0$, especially when we consider Gaussian errors. This comes partly from the
constants depending on
$\sigma$ that could completely change when $\sigma$ becomes small, and from
the bound
$$\int_{0}^{\ell_m}\exp(\sigma^2 x^2) dx\leq \int_{0}^{\ell_m}\exp(\sigma^2 \ell_m x) dx
= \frac{\exp(\sigma^2\ell_m^2) -1}{\sigma^2\ell_m}.$$
The last term is globally equivalent to $\ell_m$ when $\sigma$ tends to zero. But only the first
part $\exp(\sigma^2\ell_m^2)/(\sigma^2\ell_m^2)$
is retained for $\sigma>0$ to evaluate the rate of convergence. In a  general
setting, the dominant term for the variance term changes when $\sigma$ gets
smaller. }
\end{rem}

\section{Simulation results}
\subsection{Some examples}
\begin{figure}
$$\includegraphics[width=1\textwidth,height=5cm,]{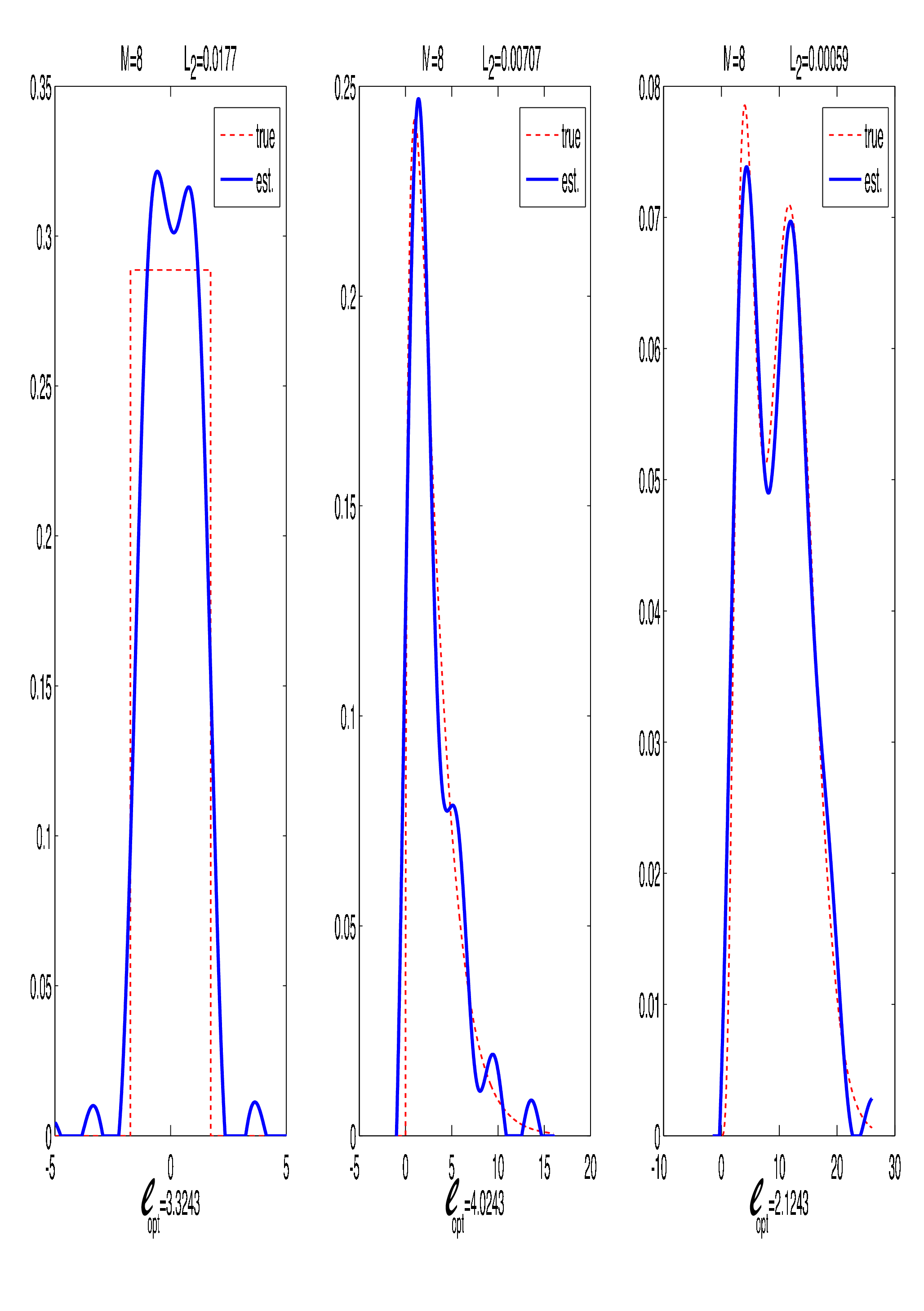}$$
\vspace{-1.3cm}\caption{Estimate and true density with $g$ uniform (left), $\chi_2(3)$ (middle) or mixed Gamma (right) with selected $\ell$ - Laplace errors -
$n=500$, $s2n$=10.}\label{ExamplesOS}
\end{figure}

\begin{figure}
$$\includegraphics[width=1\textwidth,height=5cm,]{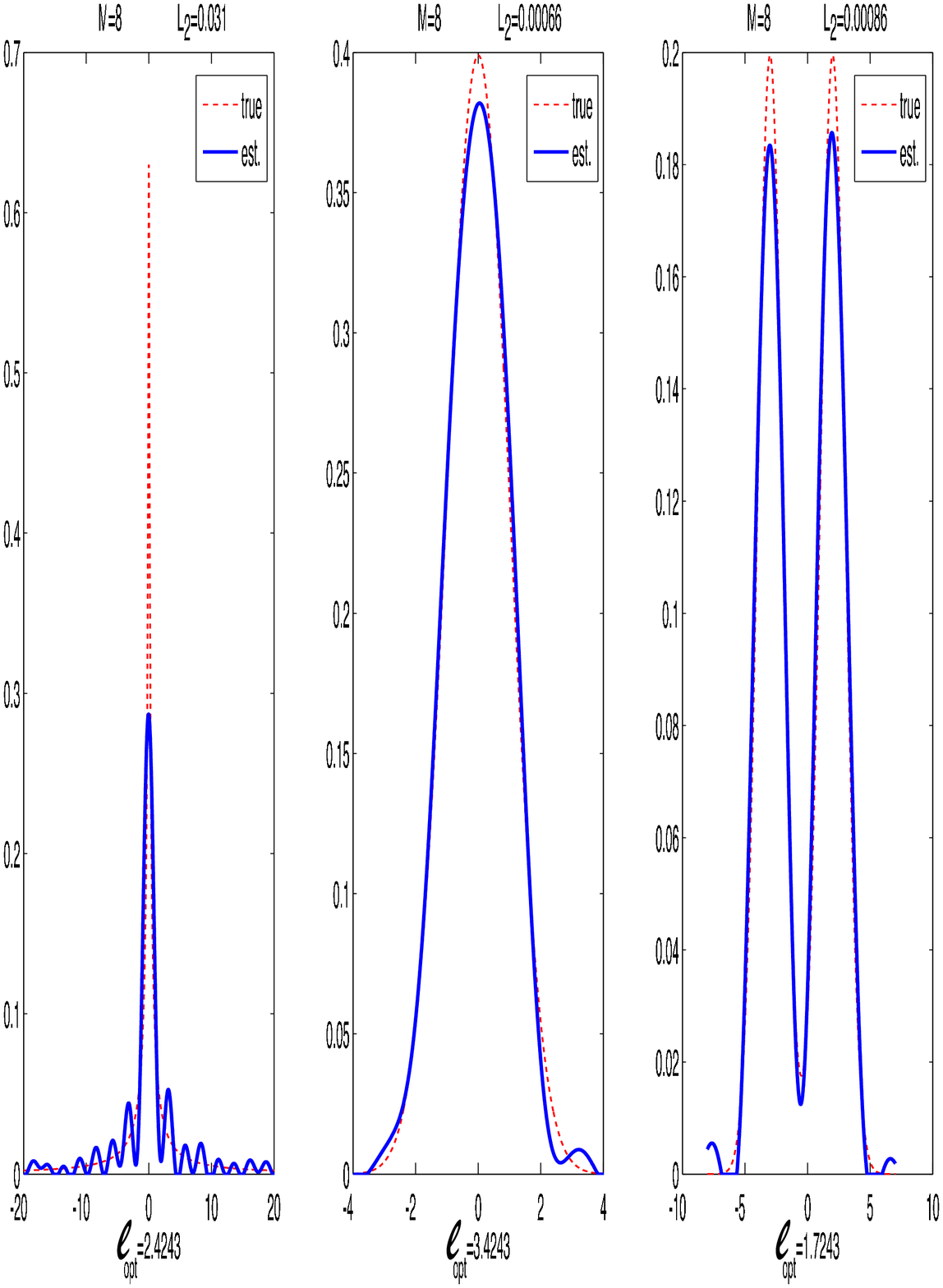}$$
 \vspace{-1.3cm}\caption{Estimate and true density with $g$ stable(1/2) (left),
Gaussian (middle) or mixed gaussian (right) with selected $\ell$
- $n=500$, $s2n$=10.}\label{ExamplesSS}
\end{figure}

Figures \ref{ExamplesOS} and \ref{ExamplesSS} illustrate the
performances of the algorithm and the quality of the estimation for
ordinary and super smooth functions $g$. Not surprisingly, the
uniform distribution or the stable 1/2 distribution are not very
well estimated, whereas the quality of the estimation for the four
other functions is very good.

Let us start a brief comparison with the results in Comte et al. ~(2005).
It is noteworthy that for the mixed gaussian
density for instance, the length selected by the algorithm with $\Delta=1/10$, corresponds to a $L_m$
which is much smaller than 1 since  $\ell_m=\pi L_m$.  Moreover, the other choices illustrate
that the algorithm takes full advantage of the more numerous
possible choices that can be done for the $\ell_m$'s. Besides, the
selected lengthes are always quite small and thus far from asymptotic.

\subsection{Mean Integrated Squared Errors}\label{misecomments} For all simulations, the MISE
is evaluated by empirical estimation over 1000 samples.
Table \ref{BasicMise} presents the MISE for the two types of errors, the
different tested densities, different $s2n$ and different sample sizes.

The first comment on Table \ref{BasicMise} concerns the importance of $\sigma$. Clearly the
MISE are smaller when there is less noise ($\sigma$ small, $s2n$ large).

The second comment is about the relative bad results for the estimation of
stable distributions, especially for stable distribution with parameter $1/4$.
If we have a look at the theoretical rate of order $(\ln(n))^{20}/n$, we easily
see that this rate tends to zero but the asymptotic is very far compared with
the considered sample sizes as it is illustrated in Section \ref{comp}.
Also note that,  in those cases, the computation of the MISE is done by using
the method (E2), which leads to larger MISE than those computed with (E1)
(two or three times (or more) larger MISE with (E2) than with (E1)), as illustrated
by the comparisons in Section \ref{compRI}.

Table \ref{BasicMise} specifies that we take $M=8$.

\subsection{Comparison of empirical and theoretical rates}
\label{comp}

The rates can be illustrated from Table \ref{BasicMise} by plotting the MISE
obtained in function of $n$. This allows to compare the empirical and the theoretical asymptotic
rates and to
evaluate the influence of the value of $\sigma^2$. It is worth emphasizing anyway that in the case
where the error is
Gaussian and $g$ super-smooth (densities \textbf{(g,l)}), the rate is directly function of $\sigma^2$.
Moreover, the rate is clearly better than logarithmic.

\begin{figure}[!ht]
$$\includegraphics[width=1\textwidth,height=18cm,]{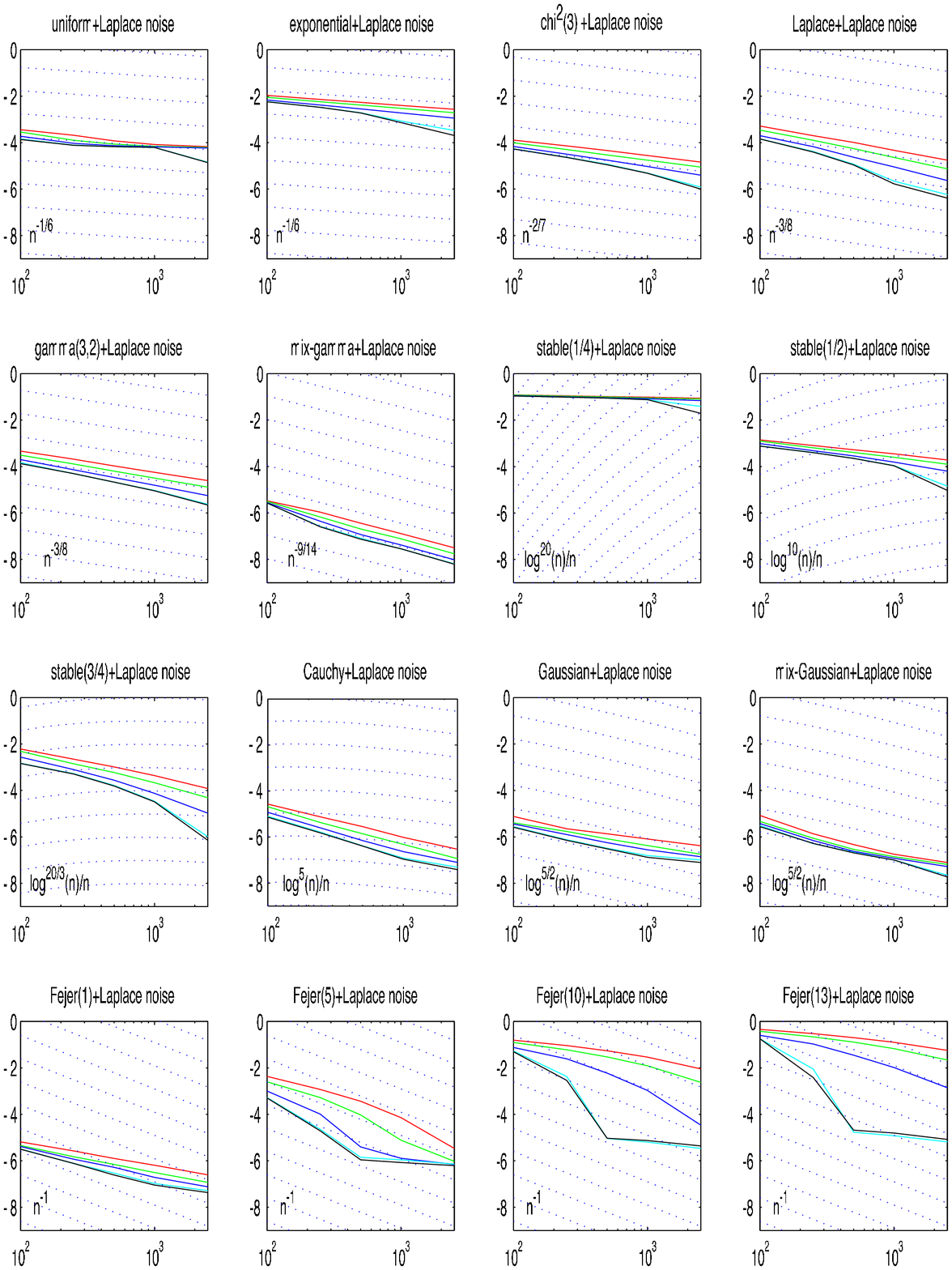}$$
\vspace{-1.3cm}\caption{Empirical MISE and theoretical asymptotical rates in logarithmic scale,
when the errors follow a Laplace distribution. From top to botton, full lines correspond
to increasing $s2n$ (2,4,10,100,1000). Dashed lines are abacuses (up to an additive constant)
for the log-theoretical rates.}\label{goodlap}
\end{figure}

\begin{figure}[!ht]
$$\includegraphics[width=1\textwidth,height=19cm,]{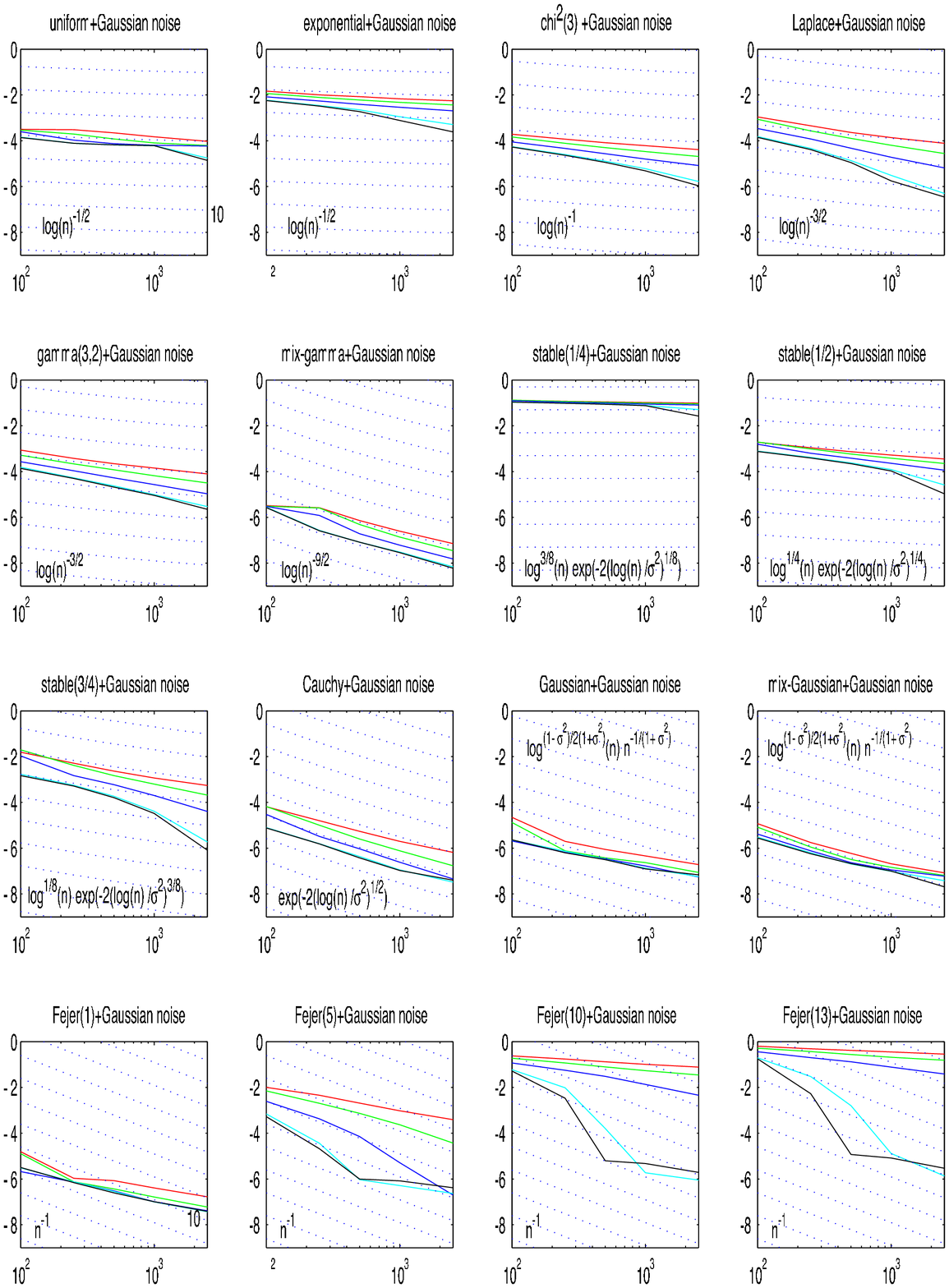}$$
\vspace{-1.3cm}\caption{Empirical MISE and theoretical  asymptotical rates in logarithmic scale,
when the errors follow a Gaussian distribution. From top to botton, full lines correspond
to increasing $s2n$ (2,4,10,100,1000). Dashed lines are abacuses (up to an additive constant)
for the log-theoretical rates.}\label{goodgauss}
\end{figure}

In order to compare the empirical MISE with the theoretical MISE, we
plot in all cases for all values of $n$ and of $s2n$, the log-MISE
in function of $\ln(n)$. In order to allow the comparison with the
theoretical rates, these log-rates are plotted with dashed lines abacuses
in function of $\ln(n)$.  Each abacus corresponds to a different value of the (unknown)
multiplicative constant in the rate.
The results are plotted in Figures
\ref{goodlap} (Laplace errors) and \ref{goodgauss} (Gaussian errors).

Consider for instance the case of Mixed Gamma distribution with Laplace errors in
 Figure \ref{goodlap}, sixth subplot. The dashed abacuses give the log of $n^{-9/14}$
(theoretical rate, see Table \ref{theorates})
up to an additive constant. The full lines give the empirical rates for $s2n=2$ to $s2n=1000$ from
top to bottom.
As $-(9/14)\ln(100)\sim -3$, one can deduce from the plot that,
since the
intercept is between -5.5 and -6, the constant is between $e^{-2.5}$ and $e^{-3}$ and the rate
of order $0.08 n^{-9/14}$ for $s2n=2$ and $0.05 n^{-9/14}$ for $s2n=1000$.

\begin{figure}
$$\includegraphics[width=0.6\textwidth,height=7cm]{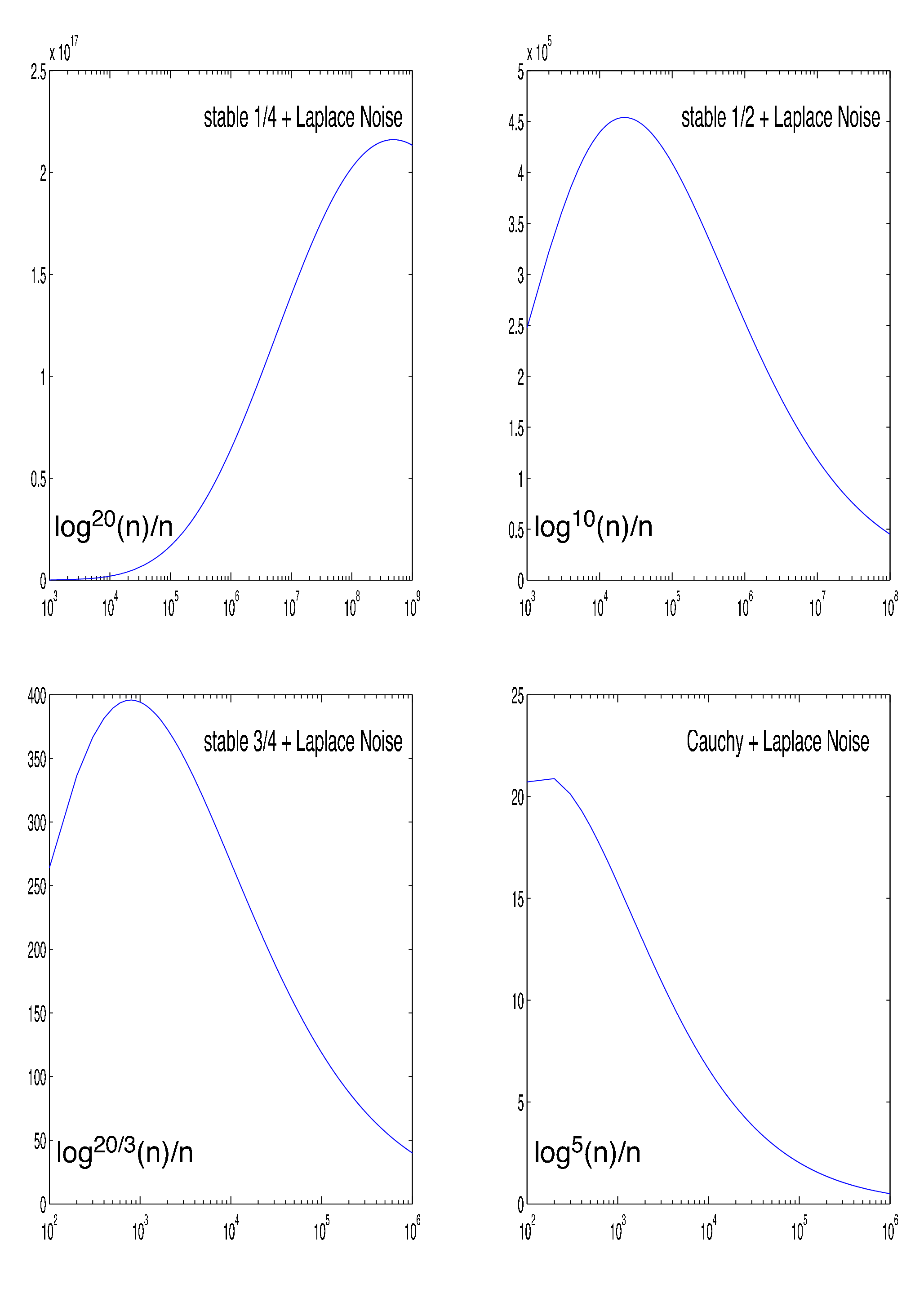}
$$ \vspace{-1.3cm}\caption{Theoretical rates in function of
$n$ for the cases where the MISE and the theoretical rates do not
seem to correspond.}\label{expltheo}
\end{figure}

\begin{figure}
$$\includegraphics[width=0.7\textwidth,height=8cm]{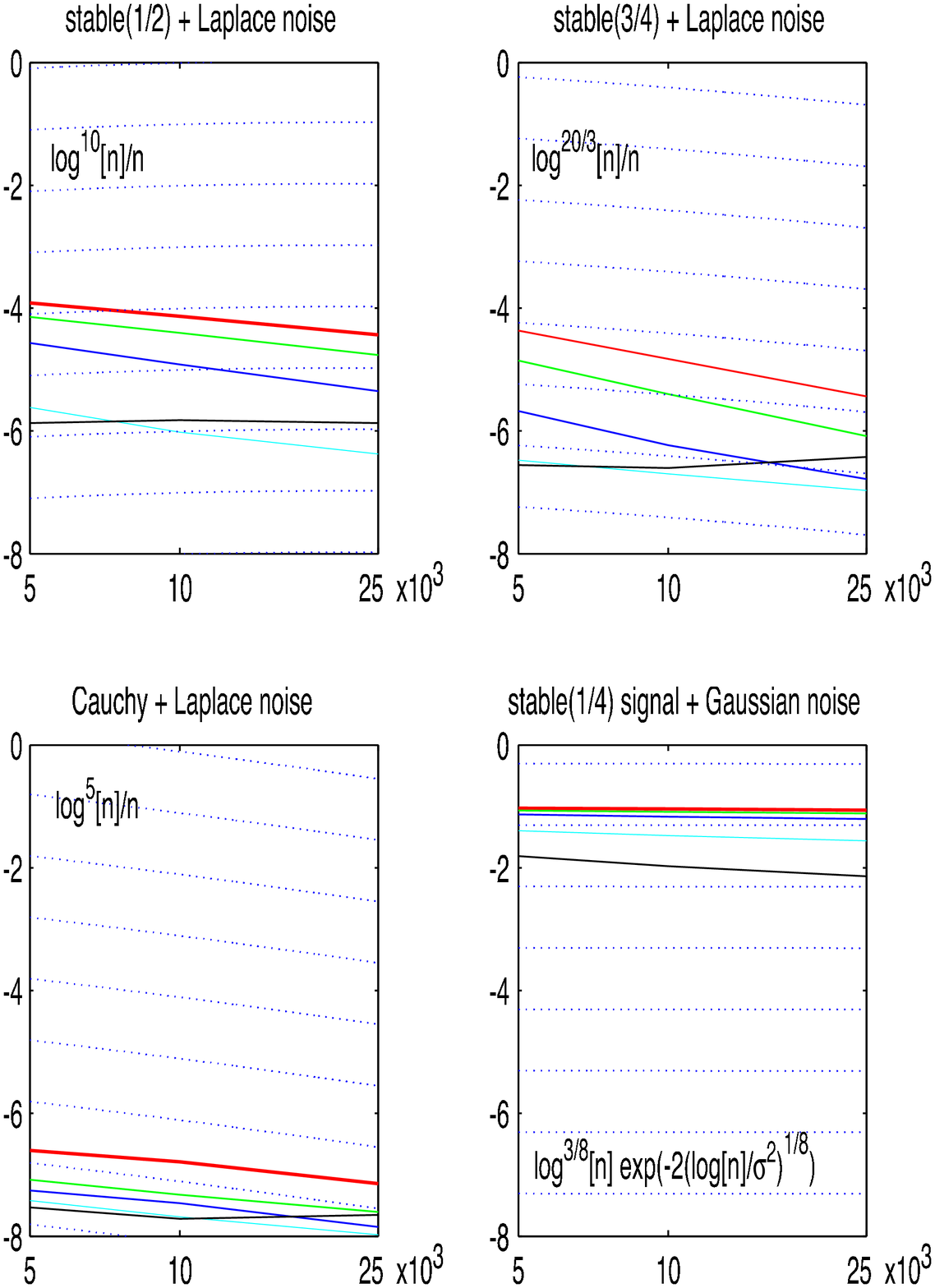}$$
\vspace{-1.3cm}\caption{Empirical MISE and theoretical rates in logarithmic scale  up to
$n=25000$ when the errors follow a Laplace distribution. The ``bad cases".}\label{badlap}
\end{figure}

We can see that most
results are in very good accordance with the theoretical predictions, but a few
results in the case of Laplace errors are less satisfactory. Figure
\ref{expltheo} explains the reason of this last fact: when we plot
the theoretical log-rates in function of $n$ in those cases, we find
out that the asymptotic that make the logarithmic part of the rate
negligible is reached for only very huge values of the sample size
$n$. It is quite positive anyway to see that in those bad cases,
our method behaves much better than what could be hoped
from the asymptotics. Figure \ref{badlap} plots these curves
including some higher values of $n$ going up to $n=25 000$, to show
how further are the asymptotics in practice.

Note that, for the rates depending on $\sigma$, we arbitrarily chose
$s2n=4$ since it was not possible to give several theoretical
curves. On the one hand, it appears from the Cauchy distribution
that even if assumption \eref{momg} is not satisfied,  the procedure
can work. On the other hand,  stable distributions show nevertheless
that a narrow pick can be  quite difficult to estimate.

\subsection{Robustness when $s2n$ is estimated}
\label{robust}
We now propose a study of the robustness
of the algorithm when
$s2n={\rm Var}(X)/\sigma^2={\rm Var}(Z)/\sigma^2-1$ is replaced by
a simple estimator (empirical variance of the observed $Z_i$'s instead
of the theoretical one). The MISE is computed with the algorithm built on a penalty with
an estimated $s2n$ with a lower bound $1/0.6$ that is about 1.67. This lower bound
is required for $s2n=2$ mainly. As we already mention it, an under-penalization
can make the MISE explode and must be avoided. We compute the ratio of the MISE
obtained with the estimated
$s2n$ over the MISE of Table \ref{BasicMise} when $s2n$ is known, and we obtain ratios equal
to one, except in the cases given in Table \ref{Ratios2n}, which remain of order one for most of
them. The empirical $s2n$ in the penalty has therefore very small influence.

\subsection{Comparison with some dependent samples}
\subsubsection{Two $\beta$-mixing examples}
In Comte et al.~(2005), most of the asymptotic properties of the adaptive
estimator $\tilde g$ are stated in the i.i.d. case, but some robustness
results are also provided. More precisely, it is shown that, when both the
$X_i$'s and the $\varepsilon_i$'s are absolutely regular, under some weak
condition on the $\beta$-mixing coefficients,  then the ${\mathbb L}_2$-risk of the
adaptive estimator $\tilde g$ has the same order as in the
independent case. The main change is the multiplicative constant in the
penalty term, which involves the sum of the $\beta$-mixing coefficients.
In other words, the adaptive procedure remains relevant for
dependent data. Here we propose to study the performances of the computed
estimator when the $X_i$'s are now $\beta$-mixing, and so are the $Z_i$'s.

This study is done  by comparing the MISE obtained respectively for the Gaussian \textbf{(k)}
and the mixed Gaussian \textbf{(l)} distributions in the independent case with the
distributions obtained in the dependent cases generated as follows.

\smallskip

$\bullet$ Construction of the dependent sequence of the $X_i$'s with
stationary standard Gaussian distribution \textbf{(k)}.

Let $(\eta_k)_{k\geq 0}$ be a sequence of i.i.d. Gaussian random variables
with mean 0 and variance $\sigma_{\eta}^2$. Let $(Y_k)_{0\leq k\leq n+1000}$ be
a sequence recursively generated by
\begin{eqnarray}
\label{dep}
Y_{k+1}= aY_k+b+\eta_{k+1}, \;\; Y_0=0,
0<a<1.\end{eqnarray} In that case, the distribution of the sequence of the $Y_k$'s
converges with exponential rate to a unique stationary distribution which is
the Gaussian distribution ${\mathcal N}(b/(1-a),\sigma_{\eta}^2/(1-a^2))$. Therefore, we take, as
an $n$-sample of $X$, the sequence $(X_1,\cdots,X_{n})=(Y_{1001}, \cdots,
Y_{n+1000})$, and we choose $b=0$, and $\sigma_{\eta}^2=1-a^2$, in \eref{dep}, so
that the resulting distribution of the $X_i$'s is the standard Gaussian
${\mathcal N}(0,1)$.
Consequently, the stationary distribution of the $X_i$'s distribution is the standard Gaussian density
\textbf{(k)}.
\\

$\bullet$ Construction of the dependent sequence of the $X_i$'s with
stationary mixed Gaussian distribution (l).

We propose here to mix two such gaussian sequences, independent from each other.
More precisely, we generate two sequences, using the method described previously.
We first generate
$Y_k^{(1)}$, $k=1,\cdots,n+1000$ with $\sigma_\eta^2=1-a^2$, $b=-3(1-a)$ and second
$Y_k^{(2)}$, $k=1,\cdots,n+1000$ with $\sigma_\eta^2=1-a^2$, $b=2(1-a)$.
Finally we generate some uniform variable on $[0,1]$, denoted by $U$ and
propose to take $X_k$ as $X_k=Y_{k+1000}^{(1)}$ if $U<0.5$ and $X_k=Y_{k+1000}^{(2)}$ else.
Clearly, the covariance between the $X_i$ and $X_{i+1}$ is divided by two
thanks to the independent additional uniform sequence standardly used for the
mixing of the distributions.
It follows that the stationary distribution of the
 $X_i$'s is the mixed Gaussian distribution \textbf{(l)}.\\

In both contexts,
we generate such sequence of $X_i$'s
for different values of $a$, $0<a<1$. Such sequences are known to be
geometrically $\beta$-mixing, with $\beta$-mixing coefficients
$(\beta_k)_{k\geq 0}$ such that $\beta_k\leq Me^{-\theta k}$, for some
constants $M$ and $\theta$. The nearer $a$ of 1, the stronger the dependency.

We study the properties of $\tilde g$, for different values of
$a$, by computing the ratio between the
resulting MISE and the MISE obtained in the independent cases \textbf{(k,l)}. The results are presented in
Table \ref{DependMise1} and Table \ref{DependMise2}.

We can see that the procedure behaves in the same way in both cases,
and that the resulting MISE ratios comparing the dependency to
independence get higher when $a$ increases and gets nearer of one.
The result remain quite good until $a=0.8$ and even $0.9$ for small $s2n$'s, if we keep in
mind that the MISE is very low in the independent case for these
two distributions.

Globally, for reasonable values of $a$ (at least
between 0 and 0.75), the dependency does not seem to bring any
additional problem.

\subsubsection{A dependent but non mixing example}

We also simulate the following dependent model. Generate $(\eta_i)_{1\leq i\leq n+1000}$
an i.i.d. Bernoulli sequence ($\eta_1=0$ or 1 with probability 1/2). Then generate
$U_{i+1}=(1/2)U_i +\eta_{i+1}$ with $U_0=0$, for $i=1, \dots, n+1000$. Take $X_k=\sqrt{3}(
U_{k+1000}-1)$  for $k=1, \dots, n$. The stationary distribution of the $U_k$'s is a uniform density
on $[0,2]$ and therefore the distribution of the $X_i$'s is the distribution \textbf{(a)}, uniform
on $[-\sqrt{3}, \sqrt{3}]$.
This model is however known to be dependent and non mixing (see
e.g. Bradley~(1986)). We experiment the estimation procedure
and we compute the ratio of the MISE for this model with the MISE in the independent case
\textbf{(a)}, for the different values of $s2n$ and sample sizes. The resulting table is not given
here because it contains essentially ones, the non ones number being at most 1.1. This may be due to
the poor quality of our estimation of the uniform distribution even in the independent context
which is then not worse in this
special dependent context. But this shows also that the procedure may be robust to some
form of dependency  quite different of the one usually met in the statistical literature.

\subsection{Comparison with Delaigle and Gijbels'(2004a)}

We propose here to compare the performances of our adaptive estimator with the
performances of the deconvolution kernel as presented in  Delaigle
and Gijbels~(2004a). This comparison  is done for densities \textbf{(e,f,k,l)}
which correspond to the densities $\#2$, $\#6$, $\#1$ and $\#3$ respectively, in  Delaigle
and Gijbels~(2004a). They
give median ISE obtained with kernel estimators by using four
different methods of bandwidth selection. The comparison is given in
Table \ref{DelGibcomp} between the median ISE computed for 1000 samples generated
with the same length and signal to noise ratio as Delaigle and
Gijbels (2004a). We compute the MISE's with direct approximation of
the integrals on the same intervals as they do, see Section
\ref{lesdensites}. We also give our corresponding means since we
think that they are more meaningful than medians. With a
multiplicative constant in the penalty smaller than the one we
chose, it may happen that medians are much better but means become
huge simply because of a few number of bad paths. The cost of such
bad paths seems therefore to have a price given by means and completely
hidden by medians.

We can see that our estimation procedure provides results of the
same quality for the ordinary smooth densities, namely for the
$\chi^2(3)$ and the Mixed Gamma densities, but that our results are
globally quite better for super-smooth densities (namely, the
Gaussian and the mixed Gaussian densities). It is noteworthy that in this case
the new penalty functions given in \eref{penlap} and \eref{pengaus} give
better MISE than the penalty functions \eref{oldpenlap} and \eref{oldpengaus} provided in Comte et
al.~(2005).

\subsection{Comparison with direct density estimation when $s2n$ is large}\label{Dalel}
We propose now to study the robustness of our procedure when $s2n$
is large, that is when the $X_i$'s are in fact almost observed. We propose to compare the non asymptotic properties of our
deconvolution estimator when $s2n=10000$,  with those, presented in a recent work by
Dalelane~(2004),
about adaptive data driven kernel estimator for density estimation, (based on the sample
$(X_1,\cdots,X_n))$. We consider here three of the
four densities considered by Dalelane~(2004), namely the normal
density \textbf{(k)}, the scale transform of the F\'ejer-de la Vall\'ee
Poussin density, the F\'ejer 5 distribution  given by \textbf{(n)} and the
$\Gamma(2,3/2)$ distribution \textbf{(d)}. The results are given in Table
\ref{DelalaneComp}. We give the MISE for Laplace errors since the
MISE for Gaussian errors are essentially the same when $s2n=10000$.

Even in these circumstances which are very unfavorable to our
estimator, we find out that our method performs very well for the
Gaussian distribution (even often better than Dalelane's~(2004)
estimator), quite well for the Gamma density where the MISE's are of
the same order, and also for the F\'ejer 5 for $n=500$ or $n=1000$.
Only the results for the F\'ejer 5 distribution when $n$ is small
($n=50, 100$) give much higher MISE's.

Therefore, it appears that our density deconvolution estimator
performs quite well despite the great number of additional numerical
approximations as compared to Dalelane's~(2004) results.

\subsection{Comparison of methods (E1) and (E2): evaluation of the MISE on
${\mathbb R}$ versus on an interval}\label{compRI}

Here, we want to compare the two methods of computation of the MISE
on an interval and on ${\mathbb R}$ as described in section
\ref{descri}, for a set of densities for which both methods are
possible: exponential, $\chi^2(3)$, Laplace, Cauchy. In those cases, we can evaluate the
bias as follows:

$$\|g-g_{m}\|^2=\frac 1{2\pi}\int_{|x|\geq \ell_m}|g^*(x)|^2dx$$
with

\noindent * for $g$ an exponential distribution \textbf{(b)},
$\int_{|x|\geq \ell_m}|g^*(x)|^2dx=2{\rm Arctan}\left(1/\ell_m\right).$

\noindent * for $g$ a normalized $\chi^2(3)$ \textbf{(c)},
$$\int_{|x|\geq \ell_m}|g^*(x)|^2dx=\sqrt{6}\left(1-\frac{2\sqrt{6}\ell_m}{\sqrt{1+(2\sqrt{6}
\ell_m)^2}}\right).$$

\noindent * for $g$ a normalized Laplace density \textbf{(d)},
$$\int_{|x|\geq \ell_m}|g^*(x)|^2dx=\sqrt{2}\left({\rm Arctan}\left(\frac{\sqrt{2}}{\ell_m}\right)
-\frac{\ell_m/\sqrt{2}}{1+ \ell_m^2/2}\right).$$

\noindent * for $g$ a Cauchy distribution \textbf{(j)},
$\int_{|x|\geq \ell_m}|g^*(x)|^2dx=e^{-2\ell_m}.$

This allows to apply method (E2) to compute the ``true" MISE on the whole real line.

It appears from Table \ref{E1E2comp} that the computation of the MISE's with method (E2)
gives results which are about two or three times greater than with method (E1), except
in the case of the exponential law where some numerical problems seem to occur when
$s2n$ becomes greater and for the $\chi_2(3)$ distribution where small samples or high levels
of noise seem to induce ratios of order 10. In the other cases, the ratio decreases when $s2n$
gets greater. The difference between the two methods of evaluation comes of course
from the oscillations of the estimate over the whole real line, even when the true
function tends to zero.

\subsection{Results when  the errors density is misspecified}
We propose here to study the non asymptotic properties of the
estimator when the error density is not correctly specified. For
both type of errors, we study the behavior of the estimator using
one type of the error density to choose the penalty when the other type of errors density
is used for the simulations of the $Z_i$'s. Table \ref{errorinerror} presents the ratio between
the resulting MISE if the errors density is not correct with the
MISE if the errors density is correct. For instance, in the first column,
the errors are Laplace but the estimator is constructed as if the
error density were Gaussian.
Some theoretical results on the effect of misspecifying the errors
distribution can be found in Meister~(2004).

\begin{figure}
$$\includegraphics[width=0.6\textwidth, height=4.5cm, ]{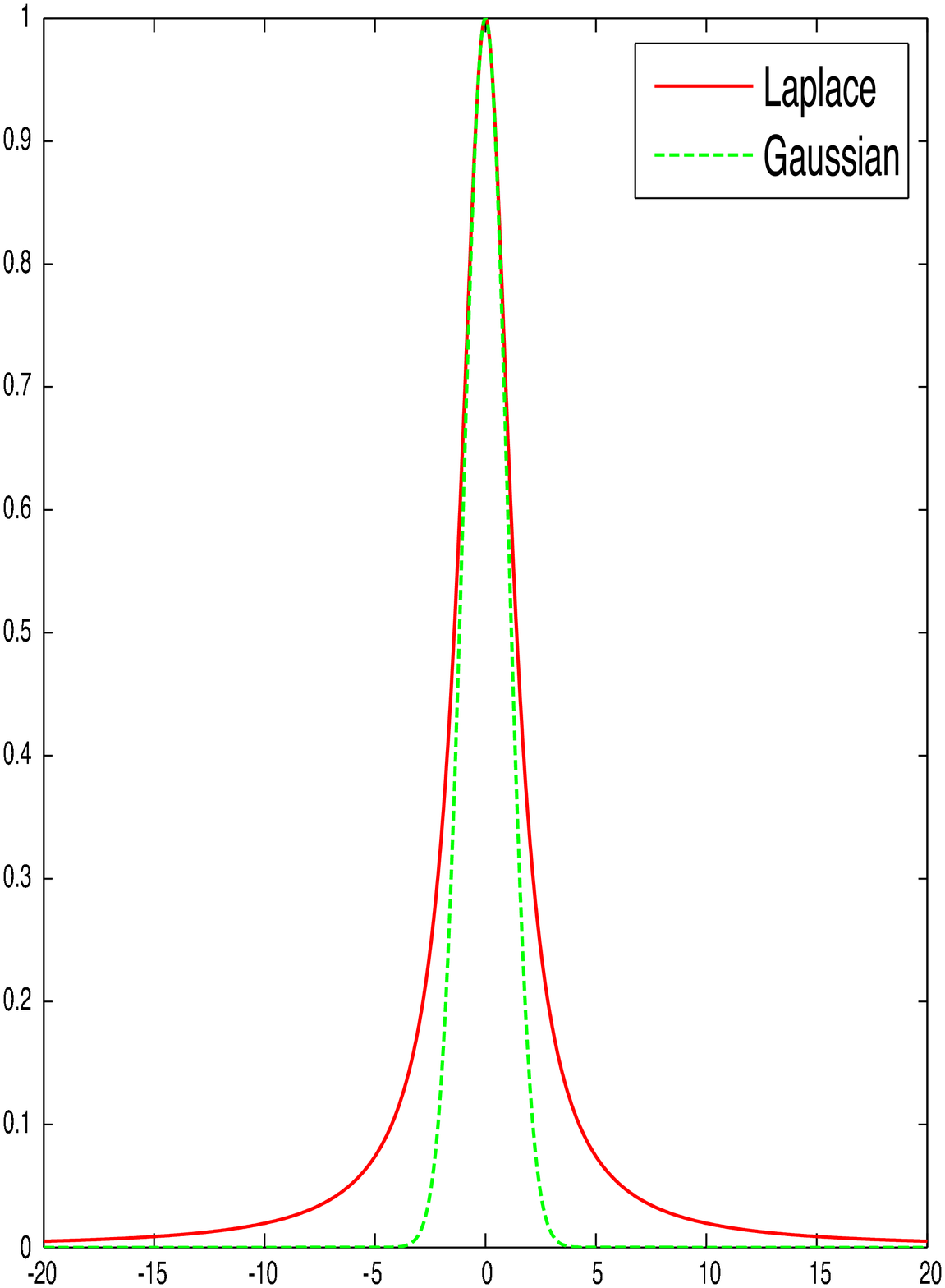}$$
\vspace{-1.3cm}\caption{Fourier transform of Laplace and Gaussian errors density.}\label{TF2}
\end{figure}

Some comments follow. As expected, since the construction uses the knowledge of the error density,
if it is misspecified, the estimator presents some bias and the MISE becomes
slightly larger. Nevertheless, this difference does not clearly appear
when $n$ is not very large. Indeed in that case, the optimal length $\ell_m$ is
small and therefore the variance term of order $\int_0^{\ell_m} \vert
f_\varepsilon^*(\sigma x)\vert^{-2}dx$ is not so quite different between the two errors.
In order to underline our comments we present in Figure \ref{TF2},
the Fourier transform of the two error densities, the Laplace and the Gaussian density.
Here, $\sigma$ is known. Globally, if we hesitate between Laplace and Gaussian errors,
Table \ref{errorinerror} seems to indicate that until $n=1000$, it is a good strategy to
always choose Gaussian errors for the estimation procedure.

We also study the behavior of our algorithm when ignoring the
noise, that is  by using our algorithm with $\sigma=0$ when $\sigma$ is not null. This amounts to consider that the $X_i$'s are observed ($Z_i=X_i$) when it
is not the case. In order to do this comparison, we simulate noisy data
($s2n=2$, 4, 10) and  run the estimation procedure as if $\sigma=0$ by putting
$s2n=10000$ in the associated penalty.
Table \ref{noerrorwitherror} presents the ratios between MISE resulting from
the procedure used with $s2n=10000$
and MISE resulting from the normal procedure which uses the knowledge of
$\sigma$ and then $s2n$.

Surprisingly, one can remark two different behaviors of the ratios on Table
\ref{noerrorwitherror}. No deterioration and even improvements for small
values of $n$. This can be explained by the fact that the penalty is smaller
when $\sigma=0$ so the algorithm can choose larger $\ell_m$ which may be of
interest for certain densities when $n$ is small. For larger values of $n$, we
clearly see an improvement to use our deconvolution algorithm against a direct
density estimation ignoring the noise.

\section{Concluding remarks}

As a conclusion, let us emphasize that we provide a complete
simulation study involving all types of possible theoretical
behaviors and rates, which are very various in the context of
density deconvolution,  depending on the type of the errors and of
the distribution to be estimated. The results are obtained with a
fast algorithm using in particular the well-known good
performances of IFFT, and are globally very satisfactory, as
compared with some other results given in the literature. The
method is very stable and reliable, even when some conditions set
by the theory are violated (as in the case of stable
distributions), and is robust to dependency in the variables. The
standard way of computing the ISE on an interval is nevertheless
proved to be more favorable than a more global method that can be
implemented here. Nevertheless the first method is the standard
one. The procedure seems also robust to a misspecification of the
error density provided that the level of the noise is well
calibrated, and is numerically stable enough to recover good
orders as compared to direct density estimation in spite of much
more (and useless in a case of direct estimation) computations.
Therefore, our global results show that the procedure works very
well, even for finite sample leading to selected lengthes very far
from the asymptotic orders.

\section*{Appendix : proof of Theorem \ref{genepenind}}

The proof essentially follows the lines of the proof of Theorem 4.1 and 4.2 in Comte et al.~(2005), and
details the role of $\Delta$. We define
$\nu_n(t)=\frac 1n\sum_{i=1}^n [u_t^*(Z_i)
-\langle t,g\rangle]$ and $B_{m,m'}(0,1)=\{t\in S_{\ell_m\vee \ell_m'}^{(n)} \;/\; \|t\|=1\}.$
Arguing as in Comte et al.~(2005), for $x>1$ we have
\begin{eqnarray*} \|\tilde g-g\|^2
\leq \left(\frac{x+1}{x-1}\right)^2 \|g-g_m^{(n)}\|^2 +
\frac{x(x+1)}{x-1}\!\!\!\sup_{t\in B_{m,\hat m}(0,1)}\!\!\!\nu_n^2(t)
+\frac{x+1}{x-1} ({\rm pen}(\ell_m)- {\rm pen}(\ell_{\hat m})).
\end{eqnarray*}
Choose some positive function $p(\ell_m,\ell_{m'})$ such that $x p(\ell_m,\ell_{m'})
\leq {\rm pen}(\ell_m) + {\rm pen}(\ell_{m'})$. Consequently, for $\kappa_x=(x+1)/(x-1)$ we have
\begin{multline}
\label{eq1}\|\tilde g-g\|^2
  \leq \kappa_x^2
  \left[\|g-g_m\|^2+\|g_m-g_m^{(n)})\|^2\right] +   x \kappa_x W_n(\ell_{\hat m})
\\
  +\kappa_x\left(x p(\ell_m,\ell_{\hat m}) +{\rm pen}(\ell_m)- {\rm pen}(\ell_{\hat m})\right)
\end{multline}
with $W_n(\ell_{m'}):=
[\sup_{t\in B_{m,m'}(0,1)} |\nu_n(t)|^2-p(\ell_m, \ell_{m'})]_+,$
and hence
\begin{equation}\label{majo2} \|\tilde g-g\|^2
  \leq \kappa_x^2 \|g-g_m\|^2+\kappa_x^2\frac{(M_2+1)\ell_m^2}{\pi^2 K_n}
  + 2\kappa_x{\rm pen}(\ell_m) + x\kappa_x\sum_{m'\in
    \mathcal{M}_n}W_n(\ell_{m'}).\end{equation} The main point of the
proof lies in studying $W_n(\ell_{m'})$,
and more precisely in finding $p(\ell_m,\ell_{m'})$ such that for a constant $K$,
\begin{equation}\label{but}
\sum_{m'\in {\mathcal M}_n} \mathbb{E}(W_n(\ell_{m'}))\leq K/(n\Delta).
\end{equation}
In that case, combining (\ref{majo2}) and (\ref{but}) we infer that, for all
$m$ in ${\mathcal M}_n$,
\begin{equation}\label{compromi}
\mathbb{E}\|g-\tilde g\|^2 \leq C_x\inf_{m\in {\mathcal M}_n}
\left[ \|g-g_m\|^2 + {\rm pen}(\ell_m)+\frac{(M_2+1)\ell_m^2}{\pi^2 K_n}\right] + x\kappa_x
\frac {K}{n\Delta},
\end{equation}
where $C_x=\kappa_x^2\vee 2\kappa_x$ suits.
It remains thus to find
$p(\ell_m,\ell_{m'})$ such that \eref{but} holds. This is done by
applying a version of Talagrand's Inequality (see Talagrand~(1996)),  to the class of
functions ${\mathcal F}= B_{m,m'}(0,1)$. If we denote by $\ell_{m^*}=\ell_m\vee
\ell_{m'}$, we get that
$$\sum_{m' \in \mathcal{M}_n}\mathbb{E}(W_n(\ell_{m'}))\leq K\sum_{m' \in
  \mathcal{M}_n}[I(\ell_{m^*})+II(\ell_{m^*})],$$
where $I(\ell_{m^*})$ and $II(\ell_{m^*})$ are defined
by
\begin{eqnarray*}
&&I(\ell_{m^*})=\frac{\lambda_2\ell_{m^*}^{2\gamma+
(1/2-\delta/2)\wedge(1-\delta)}\exp\{2\mu\sigma^\delta \ell_{m^*}^{\delta}\}}{n}
\exp\{{-K_1\xi^2(\lambda_1/\lambda_2)\ell_{m^*}^{(1/2-\delta/2)_+}}\},\\
&& II(\ell_{m^*})= \frac{\lambda_1
\ell_{m^*}^{2\gamma+1-\delta}e^{2\mu\sigma^\delta(\ell_{m^*})^{\delta }}}{n^2}
\exp\left\{{-(K_1\xi C(\xi) \sqrt{n}/\sqrt{2}}\right\},
\end{eqnarray*}
with for $\ell_m\geq \ell_0$, $$\lambda_2=\left\{ \begin{array}{ll} 1 & \mbox{if } \delta>1 \\
\lambda_1^{1/2} (\ell_0^{-2}+ \sigma^2)^{\gamma/2} \|f_{\varepsilon}^*\| \kappa_0^{-1}
(2\pi)^{-1/2} & \mbox{ if } \delta\leq 1.\end{array}\right.$$

\noindent {\bf 1) Study of $\sum_{m\in \mathcal{M}_n}II(\ell_{m^*})$}.\\
If we denote by $\Gamma(\ell_m)=
\ell_m^{2\gamma+1-\delta}\exp\{2\mu\sigma^\delta \ell_m^{\delta}\}$ then
\begin{eqnarray*}
\sum_{m\in \mathcal{M}_n}II(\ell_{m^*})\leq C(\lambda_1)\vert \mathcal{M}_n\vert \exp\left\{-
(K_1\xi C(\xi)\sqrt{n})/\sqrt{2}\right\}\Gamma(\ell_{m_n})/n^2.
\end{eqnarray*}
Consequently, as soon as $\Gamma(\ell_{m_n})/n$ is bounded (we only consider
$m_n$ such that pen$(\ell_{m_n})$ is bounded), then $\sum_{m\in  \mathcal{M}_n}II(\ell_{m^*})\leq C/n$

\noindent {\bf 2) Study of $\sum_{m\in \mathcal{M}_n}I(\ell_{m^*})$.}\\
Denote by
$\psi=2\gamma+ (1/2-\delta/2)\wedge(1-\delta)$, $\omega=(1/2-\delta/2)_+$,
$K'=K_1\lambda_1/\lambda_2$, then for $a,b\geq 1$, we infer that

\begin{eqnarray}\nonumber
(a\vee b)^{\psi}e^{2\mu\sigma^\delta (a \vee b)^{\delta}}e^{-K'\xi^2(a\vee b)^{\omega}}
&\leq&
(a^{\psi}e^{2\mu\sigma^\delta a^{\delta}}+b^{\psi} e^{2\mu\sigma^\delta
b^{\delta}})e^{-(K'\xi^2/2)(a^{\omega} +
b^{\omega})}\\&\leq&   a^{\psi}e^{2\mu\sigma^\delta
a^{\delta}}e^{-(K'\xi^2/2)a^{\omega}} e^{-(K'\xi^2/2)
b^{\omega}}+b^{\psi} e^{2\mu\sigma^\delta b^{\delta}} e^{-(K'\xi^2/2)
b^{\omega}}.\label{eqmax}\end{eqnarray}
Consequently, if we denote by $\tilde \Gamma$ the quantity
$\tilde \Gamma(\ell_{m^*})=\ell_{m^*}^{2\gamma+
(1/2-\delta/2)\wedge(1-\delta)}\exp\{2\mu\sigma^\delta \ell_{m^*}^{\delta}\}$
then
\begin{eqnarray}
\sum_{m' \in \mathcal{M}_n}I(\ell{m^*})
&\leq& C_1(\lambda_2)\frac{\tilde \Gamma(m)}{n}\exp\{-(K'\xi^2/2)\ell_m^{(1/2-\delta/2)}\}\sum_{m' \in
  \mathcal{M}_n} \exp\{-(K'\xi^2/2)\ell_{m'}^{(1/2-\delta/2)}\}\nonumber\\
&&+ C_1(\lambda_2)\sum_{m' \in \mathcal{M}_n}\frac{\tilde\Gamma(\ell_{m'})}{n}\exp\{-(K'\xi^2)
\ell_{m'}^{(1/2-\delta/2)}\}\label{I}.
\end{eqnarray}

\paragraph{\textbf{a) Case $0\leq \delta < 1/3$}} In that case, since
$\delta< (1/2-\delta/2)_+$, the choice $\xi^2=1$ ensures that
$\tilde \Gamma(\ell_m)\exp\{-(K'\xi^2/2)(\ell_m)^{(1/2-\delta/2)}\}$ is bounded and
thus the first term in \eref{I} is bounded by
$$\frac{C}{n\Delta}\int_0^\infty \exp\{-(K'\xi^2)x^{(1/2-\delta/2)}\}dx\leq
\tilde C/(n\Delta).$$
In the same way, $\sum_{m' \in
   \mathcal{M}_n}\tilde\Gamma(\ell_{m'})\exp\{-(K'\xi^2)\ell_{m'}^{(1/2-\delta/2)}\}/n$ is bounded by
\begin{multline*}
\frac{C}{n\Delta}\int_0^\infty (x+1)^{2\gamma+ (1/2-\delta/2)\wedge(1-\delta)}
\exp\{2\mu\sigma^\delta( (x+1))^{\delta}\}\exp\{-(K'\xi^2)x^{(1/2-\delta/2)}\}dx
\leq \tilde{\tilde C}/(n\Delta).
\end{multline*}
It follows that
$\sum_{m'\in {\mathcal M}_n}I(\ell_{m^*})\leq C/(n\Delta).$
Consequently, \eref{but}
holds if we choose $\mbox{pen}(\ell_m)=2x(1+2\xi^2) \lambda_1
\ell_m^{2\gamma+1-\delta}\exp\{2\mu\sigma^\delta \ell_m^{\delta}\}/n.
$
\paragraph{\textbf{b) Case $\delta= 1/3$}} In that case, bearing in mind
Inequality \eref{eqmax} we choose
$\xi^2$
such that $2\mu\sigma^\delta \ell_{m^*}^{\delta} -
(K'\xi^2/2)\ell_{m^*}^{\delta}= -2\mu\sigma^\delta \ell_{m^*}^{\delta}$
that is
$\xi^2=(4\mu\sigma^\delta\lambda_2)/(K_1\lambda_1).$
By the same arguments as for the case $0\leq \delta<1/3$, this choice ensures that
$\sum_{m'\in {\mathcal M}_n}I(\ell_{m^*}) \leq C/(n\Delta)$, and consequently \eref{but}
holds.
The result follows by taking
$p(\ell_m,\ell_{m'})=2(1+2\xi^2)\lambda_1 \ell_{m^*}^{2\gamma+1-\delta}\exp(2\mu\sigma^\delta
\ell_{m^*}^{\delta})/ n,$
and
$\mbox{pen}(\ell_m )=2x(1+2\xi^2)\lambda_1 \ell_m^{2\gamma+1-\delta}\exp(2\mu\sigma^\delta
\ell_m^{\delta})/ n.$

\paragraph{\textbf{c) Case $\delta> 1/3$}} In that case, $\delta > (1/2-\delta/2)_+$.
Bearing in mind Inequality \eref{eqmax} we choose $\xi^2=\xi^2(\ell_m,\ell_{m'})$
such that $2\mu\sigma^\delta\ell_{m^*}^{\delta} -
(K'\xi^2/2)\ell_{m^*}^{\omega}= -2\mu\sigma^\delta \ell_{m^*}^{\delta}$
that is
$$\xi^2=\xi^2(\ell_m,\ell_{m'})=(4\mu\sigma^\delta\lambda_2)/(K_1\lambda_1)
\ell_{m^*}^{\delta-\omega}.$$
This choice ensures that $\sum_{m'\in {\mathcal M}_n}I(\ell_{m^*}) \leq
C/(n\Delta)$, and consequently \eref{but}
holds and \eref{resu1} follows if
$p(\ell_m,\ell_{m'})=2(1+2\xi^2(\ell_m,\ell_{m'}))\lambda_1 \ell_{m^*}^{2\gamma+1-\delta}
\exp(2\mu\sigma^\delta \ell_{m^*}^{\delta})/ n,$
and
$\mbox{pen}(\ell_m)=2x(1+2\xi^2(\ell_m,\ell_m))\lambda_1 \ell_m^{2\gamma+1-\delta}
\exp(2\mu\sigma^\delta \ell_m^{\delta})/ n.$
\hfill $\Box$

\renewcommand{\baselinestretch}{1}

\newpage

\section*{Tables and Figures}

\begin{center}
\begin{table}
\begin{tabular}{||c||c|c||c||}\hline\hline
$g$ & $\varepsilon \sim $ Laplace & $\varepsilon \sim$ Gauss  & $\varepsilon=0$ \\
 & ($\gamma=2, \delta=0$)& ($\gamma=0, \delta=2$)&   \\ \hline\hline
\textbf{(a,b)} & & & \\
Uniform, Exponential & $n^{-1/6}$ & $[\ln(n)]^{-1/2}$ & $n^{-1/2}$ \\
$s=1/2$, $r=0$ & & & \\ \hline
\textbf{(c)}   & &  & \\
 $\chi^2(3)$,  & $\displaystyle n^{-2/7}$ & $[\ln(n)]^{-1}$ & $n^{-2/3}$ \\
$s=1, r=0$ & &  & \\ \hline
\textbf{(d,e)} & & & \\
Laplace, $\Gamma(2,3/2)$  & $\displaystyle n^{-3/8}$ & $[\ln(n)]^{-3/2}$ & $n^{-3/4}$ \\
$s=3/2, r=0$ & & & \\\hline
\textbf{(f)} & &  & \\
Mixed Gamma &  $n^{-9/14}$ & $[\ln(n)]^{-4.5}$ & $n^{-9/10}$ \\
$s=9/2, r=0$ & &  & \\ \hline
\textbf{(g)} Stable 1/4 &  & & \\
$s=-3/8$, $r=1/4$, & $\displaystyle \frac{[\ln(n)]^{20}}n$ & $[\ln(n)]^{3/8}\exp\left(-2\left(\frac{\ln(n)}{\sigma^2}\right)^{1/8}\right)$
 & $\frac{\ln^4(n)}n$ \\
 $b=1$ & & & \\ \hline
\textbf{(h)} Stable 1/2 &  & & \\
$s=-1/4$, $r=1/2$,& $\displaystyle \frac{[\ln(n)]^{10}}n$ & $[\ln(n)]^{1/4}\exp\left(-2\left(\frac{\ln(n)}{\sigma^2}\right)^{1/4}\right)$
& $\frac{\ln^2(n)}n$ \\
 $b=1$ & & & \\ \hline
\textbf{(i)} Stable 3/4 &  & & \\
$s=-1/8$, $r=3/4$, & $\displaystyle \frac{[\ln(n)]^{20/3}}n$ & $[\ln(n)]^{1/8}\exp\left(-2\left(\frac{\ln(n)}{\sigma^2}\right)^{3/8}\right)$
& $\frac{\ln^{4/3}(n)}n$ \\
 $b=1$ & & & \\ \hline
\textbf{(j)} && &  \\
Cauchy, $r=1$, & $\displaystyle \frac{[\ln(n)]^5}n$ & $\exp\left( -2\sqrt{\frac{\ln(n)}{\sigma^2}}\right)$ & $\frac{\ln(n)}n $\\
$ s=0, b=1$ & & & \\ \hline
\textbf{(k,l)}  Gauss, & & & \\
 Mixed Gauss, $r=2$,  &  $\displaystyle \frac{[\ln(n)]^{5/2}}{n}$ &
$\displaystyle (\ln(n))^{-\frac 12\frac{\sigma^2 -1}{\sigma^2+1}} \left(\frac{1}{n}\right)^{1/(1+\sigma^2)}$ & $\frac{\sqrt{\ln(n)}}n$  \\
$s=1/4, b=1/2$ & &  & \\\hline
\textbf{(m,n,o)} & & &\\
F\'ejer-DVP & $n^{-1}$ & $n^{-1}$& $n^{-1}$\\
no bias  & & & \\ \hline \hline
\end{tabular}
\caption{Theoretical orders of the rates of the adaptive estimator as deduced from Table
\ref{rates} and formulae (\ref{vitextra}) and
(\ref{vitforstable}) when $\sigma>0$ and (last column) when $\sigma=0$.}\label{theorates}
\end{table}
\end{center}

\begin{center}
\begin{table}
{\small
\begin{tabular}{||c|r||c|c||c|c||c|c||c|c||c|c||}\hline\hline
\multicolumn{2}{||c||}{$\times10^{-2}$}  & \multicolumn{2}{c||}{$n=100$}& \multicolumn{2}{c||}{$n=250$}&
\multicolumn{2}{c||}{$n=500$}& \multicolumn{2}{c||}{$n=1000$}& \multicolumn{2}{c||}{$n=2500$}\\ \hline \hline
$g$ & $s2n$ & Lap. & Gaus. & Lap. & Gaus. & Lap. & Gaus. & Lap. & Gaus. & Lap. & Gaus. \\ \hline\hline
 \multirow{5}{.5cm}{\rot{Uniform}}
 &    2 & 3.55 & 3.42 & 2.62 & 3.13 & 2.07 & 2.69 & 1.75 & 2.25 &  1.58 & 1.87\\
 &    4 & 3.06 & 3.08 & 2.13 & 2.49 & 1.78 & 2.02 & 1.61 &  1.7 &   1.5 & 1.52\\
 &   10 & 2.54 & 2.83 & 1.85 & 1.94 & 1.63 & 1.65 & 1.54 & 1.53 &  1.46 & 1.47\\
 &  100 & 2.25 & 2.27 &  1.7 &  1.7 & 1.56 & 1.56 &  1.5 &  1.5 & 0.815 & 1.07\\
 & 1000 & 2.21 & 2.22 & 1.68 &  1.7 & 1.56 & 1.55 &  1.5 & 1.49 & 0.785 & 0.79\\
\hline\hline
\multirow{5}{.5cm}{\rot{Expon. }}
 &    2 & 14.2 & 16.1 & 11.9 & 13.8 & 10.5 & 12.6 & 9.11 & 11.6 & 7.75 & 10.5\\
 &    4 &   13 & 14.6 & 10.7 & 12.3 & 9.25 & 10.9 & 8.08 & 9.82 & 6.69 & 8.74\\
 &   10 & 11.6 & 12.6 &  9.3 & 10.4 & 7.89 & 9.06 & 6.57 &  7.9 & 5.27 & 6.77\\
 &  100 & 10.8 & 10.9 & 8.45 & 8.61 & 6.66 & 7.03 & 4.66 & 5.23 & 3.11 & 3.74\\
 & 1000 & 10.7 & 10.8 & 8.37 &  8.4 & 6.55 & 6.58 & 4.37 & 4.49 & 2.53 & 2.71\\
\hline\hline
\multirow{5}{.5cm}{\rot{Chi2(3)}}
 &    2 & 2.15 & 2.51 &  1.64 &     2 &  1.33 &  1.73 &  1.06 &  1.49 & 0.811 &  1.25\\
 &    4 & 1.88 & 2.22 &  1.39 &  1.67 &   1.1 &  1.38 &  0.88 &  1.15 & 0.648 & 0.923\\
 &   10 & 1.62 &  1.8 &  1.14 &  1.33 &  0.88 &  1.05 & 0.667 & 0.829 & 0.457 & 0.624\\
 &  100 & 1.45 & 1.47 &     1 &  1.03 & 0.735 & 0.758 & 0.502 & 0.547 & 0.273 & 0.315\\
 & 1000 & 1.43 & 1.44 & 0.995 & 0.995 & 0.723 & 0.726 & 0.499 & 0.499 & 0.253 & 0.259\\
\hline\hline
\multirow{5}{.5cm}{\rot{Laplace}}
 &    2 & 3.87 & 5.19 &  2.6 & 3.52 &  1.92 &  2.74 &   1.4 &  2.17 & 0.921 &   1.7\\
 &    4 & 3.24 & 4.77 & 2.09 & 2.84 &  1.48 &  2.09 &  1.01 &  1.53 &  0.63 &  1.07\\
 &   10 & 2.61 & 3.25 & 1.61 & 2.03 &  1.03 &  1.36 & 0.677 & 0.916 &  0.39 & 0.577\\
 &  100 & 2.24 & 2.33 &  1.3 & 1.36 & 0.753 & 0.798 & 0.375 & 0.422 & 0.213 & 0.199\\
 & 1000 & 2.23 & 2.22 & 1.28 & 1.29 & 0.731 & 0.733 & 0.329 & 0.339 & 0.182 & 0.171\\
\hline\hline
\multirow{5}{.5cm}{\rot{Gamma}}
 &    2 & 3.86 & 4.83 & 2.64 & 3.54 &  1.97 &  2.73 &  1.49 &  2.21 &  1.04 &  1.72\\
 &    4 & 3.17 & 3.96 & 2.12 & 2.65 &  1.55 &  2.03 &  1.16 &  1.56 & 0.767 &  1.14\\
 &   10 & 2.59 & 2.96 & 1.66 & 1.95 &  1.18 &  1.43 & 0.851 &  1.07 & 0.534 & 0.712\\
 &  100 & 2.27 & 2.31 & 1.42 & 1.45 & 0.978 &  1.01 & 0.674 & 0.692 & 0.374 & 0.408\\
 & 1000 &  2.2 & 2.22 &  1.4 & 1.42 & 0.974 & 0.968 & 0.663 & 0.661 & 0.359 & 0.361\\
\hline\hline
\multirow{5}{.5cm}{\rot{Mix.Gamma}}
 &    2 & 0.465 &  0.47 & 0.277 & 0.362 &  0.172 &  0.241 &  0.109 &  0.144 & 0.0601 & 0.0838\\
 &    4 & 0.432 & 0.428 & 0.237 & 0.352 &  0.135 &  0.206 &  0.086 &  0.112 & 0.0453 & 0.0605\\
 &   10 & 0.396 & 0.423 & 0.196 & 0.279 &  0.106 &  0.135 & 0.0664 & 0.0773 &  0.035 & 0.0427\\
 &  100 & 0.368 & 0.386 & 0.159 & 0.163 & 0.0897 &  0.091 & 0.0556 & 0.0573 & 0.0292 & 0.0299\\
 & 1000 & 0.375 & 0.368 & 0.154 & 0.158 & 0.0867 & 0.0913 & 0.0552 & 0.0557 & 0.0288 & 0.0281\\
\hline\hline
\multirow{5}{.5cm}{\rot{Stable 1/4}}
&    2 & 40.1 & 41.3 & 38.4 & 39.4 & 37.2 & 38.3 & 36.1 & 37.4 & 34.5 & 36.4\\
 &    4 & 39.7 & 41.5 & 37.7 & 39.1 & 36.4 & 37.6 & 35.1 & 36.5 & 33.4 & 35.2\\
 &   10 & 38.9 & 40.6 & 36.9 & 37.7 & 35.4 & 36.3 & 33.7 &   35 & 31.2 & 33.3\\
 &  100 & 38.2 & 38.4 & 36.3 & 36.4 & 34.6 & 34.9 & 32.3 & 32.9 & 24.4 & 27.5\\
 & 1000 & 38.1 & 38.2 & 36.3 & 36.3 & 34.6 & 34.6 & 32.2 & 32.3 & 18.1 & 20.7\\
\hline\hline
\multirow{5}{.5cm}{\rot{Stable 1/2}}
 &    2 & 5.84 & 6.71 & 4.58 & 5.34 &  3.8 & 4.42 & 3.17 &  3.8 &  2.45 &  3.19\\
 &    4 & 5.59 & 6.78 & 4.13 & 4.96 & 3.38 & 3.96 & 2.74 &  3.3 &  2.03 &  2.62\\
 &   10 & 5.02 &  6.1 & 3.66 & 4.14 & 2.91 &  3.3 & 2.24 & 2.64 &  1.51 &  1.94\\
 &  100 & 4.55 & 4.64 & 3.37 & 3.42 & 2.62 & 2.68 & 1.91 & 1.99 & 0.791 &  1.03\\
 & 1000 & 4.51 & 4.52 & 3.34 & 3.34 &  2.6 &  2.6 &  1.9 &  1.9 & 0.661 & 0.704\\
\hline\hline
\multirow{5}{.5cm}{\rot{Stable 3/4}}
 &    2 & 10.9 & 16.8 & 7.04 & 10.2 &    5 & 7.32 & 3.46 & 5.34 &     2 &  3.78\\
 &    4 & 10.1 & 17.7 & 5.81 & 9.49 & 3.96 & 5.98 & 2.55 & 4.06 &  1.34 &  2.47\\
 &   10 &    8 & 13.7 & 4.51 & 6.11 & 2.85 & 3.93 & 1.63 & 2.45 & 0.691 &  1.23\\
 &  100 &  6.2 & 6.64 & 3.79 & 3.96 & 2.29 & 2.42 & 1.16 & 1.26 & 0.262 & 0.337\\
 & 1000 & 6.08 & 6.17 & 3.72 & 3.77 & 2.25 & 2.27 & 1.13 & 1.14 & 0.225 & 0.235\\
\hline\hline
 \end{tabular}}
 \end{table}
\end{center}

\begin{center}
\begin{table}
{\small
\begin{tabular}{||c|r||c|c||c|c||c|c||c|c||c|c||}\hline\hline
\multicolumn{2}{||c||}{$\times10^{-2}$}  & \multicolumn{2}{c||}{$n=100$}& \multicolumn{2}{c||}{$n=250$}&
\multicolumn{2}{c||}{$n=500$}& \multicolumn{2}{c||}{$n=1000$}& \multicolumn{2}{c||}{$n=2500$}\\ \hline \hline
$g$ & $s2n$ & Lap. & Gaus. & Lap. & Gaus. & Lap. & Gaus. & Lap. & Gaus. & Lap. & Gaus. \\ \hline\hline
\multirow{5}{.5cm}{\rot{Cauchy}}
&    2 &   1.2 &  1.62 & 0.683 & 0.935 & 0.449 & 0.606 & 0.294 & 0.382 &  0.185 &  0.243\\
 &    4 &  1.04 &  1.64 &  0.52 & 0.714 & 0.319 & 0.397 & 0.208 & 0.238 &  0.118 &  0.128\\
 &   10 & 0.816 &  1.18 & 0.411 & 0.458 & 0.238 & 0.265 & 0.151 & 0.149 & 0.0947 & 0.0751\\
 &  100 & 0.695 & 0.701 & 0.335 & 0.338 & 0.192 & 0.189 &  0.11 & 0.108 & 0.0736 & 0.0624\\
 & 1000 & 0.671 &  0.69 & 0.321 & 0.338 & 0.186 &  0.18 & 0.107 & 0.104 & 0.0674 &  0.067\\
\hline\hline
\multirow{5}{.5cm}{\rot{Gauss.}}
 &    2 & 0.928 &  1.09 & 0.537 & 0.538 & 0.416 & 0.397 & 0.314 & 0.281 &   0.23 &  0.194\\
 &    4 & 0.649 & 0.838 & 0.415 & 0.312 & 0.305 & 0.225 & 0.226 & 0.173 &  0.149 &  0.108\\
 &   10 & 0.609 &  0.48 &  0.37 &  0.28 & 0.248 &   0.2 & 0.178 & 0.146 &  0.128 & 0.0857\\
 &  100 & 0.522 & 0.501 & 0.283 &  0.28 &  0.19 & 0.187 & 0.133 & 0.122 &  0.101 & 0.0827\\
 & 1000 & 0.489 & 0.488 & 0.262 &  0.26 & 0.179 & 0.181 & 0.123 & 0.121 & 0.0911 & 0.0848\\
\hline\hline
\multirow{5}{.5cm}{\rot{Mix. Gauss.}}
 &    2 & 0.727 &  0.82 & 0.337 & 0.378 &   0.2 & 0.222 &  0.132 &  0.142 & 0.0892 & 0.0915\\
 &    4 & 0.562 & 0.668 & 0.267 & 0.297 & 0.167 &  0.17 &  0.115 &  0.115 & 0.0788 &   0.08\\
 &   10 & 0.498 & 0.529 & 0.242 & 0.244 & 0.151 & 0.147 &  0.107 &  0.103 & 0.0678 & 0.0762\\
 &  100 & 0.471 & 0.459 & 0.213 & 0.225 & 0.141 & 0.139 &    0.1 & 0.0983 & 0.0516 & 0.0553\\
 & 1000 & 0.453 & 0.457 & 0.216 & 0.224 & 0.141 & 0.142 & 0.0991 & 0.0979 & 0.0491 &   0.05\\
\hline\hline
\multirow{5}{.5cm}{\rot{F\'ejer 1}}
 &    2 & 0.884 &  1.02 & 0.531 & 0.426 & 0.372 & 0.393 & 0.276 & 0.262 &  0.191 &  0.181\\
 &    4 & 0.655 & 0.813 & 0.388 & 0.285 & 0.271 & 0.215 & 0.196 & 0.146 &   0.12 &  0.094\\
 &   10 & 0.616 & 0.465 & 0.341 & 0.281 &  0.23 & 0.185 & 0.147 & 0.117 & 0.0962 & 0.0736\\
 &  100 & 0.522 & 0.504 & 0.274 & 0.269 &  0.17 & 0.172 & 0.112 & 0.104 & 0.0766 & 0.0667\\
 & 1000 & 0.516 & 0.514 & 0.262 & 0.262 &  0.16 & 0.164 & 0.102 & 0.104 & 0.0684 & 0.0676\\
\hline\hline
\multirow{5}{.5cm}{\rot{F\'ejer 5}}
 &    2 & 9.43 & 13.7 &  5.53 &  9.33 &  3.29 &  6.85 &  1.66 &  5.04 & 0.557 &  3.41\\
 &    4 &  7.5 & 11.6 &  3.81 &  6.66 &  1.87 &  4.28 & 0.672 &  2.62 &  0.32 &  1.24\\
 &   10 & 5.16 & 7.52 &  1.98 &  3.56 & 0.556 &  1.66 & 0.361 & 0.543 & 0.273 &  0.16\\
 &  100 & 3.92 & 4.24 &  1.03 &  1.22 & 0.363 & 0.311 &  0.32 & 0.229 & 0.243 & 0.157\\
 & 1000 & 3.83 &  3.9 & 0.969 & 0.982 & 0.319 & 0.319 & 0.273 & 0.265 & 0.221 & 0.188\\
\hline\hline
\multirow{5}{.5cm}{\rot{F\'ejer 10}}
 &    2 & 44.5 & 53.5 & 35.2 & 46.1 &  27.9 &  41.4 &  21.3 &  37.3 &    13 &  32.9\\
 &    4 & 40.1 & 48.5 & 29.3 & 38.9 &  21.7 &  33.3 &  14.6 &  28.3 &  7.37 &  22.9\\
 &   10 & 32.5 & 39.5 & 20.1 & 28.5 &  10.8 &  21.5 &  5.21 &  15.6 &  1.24 &  9.49\\
 &  100 & 27.7 & 29.1 & 9.47 & 13.3 & 0.854 &  2.37 & 0.736 & 0.418 &  0.52 &   0.3\\
 & 1000 & 27.4 & 27.7 & 8.12 & 8.74 & 0.815 & 0.709 & 0.695 & 0.568 & 0.522 & 0.371\\
\hline\hline
\multirow{5}{.5cm}{\rot{F\'ejer 13}}
 &    2 & 70.6 & 81.7 & 59.1 & 73.1 & 49.5 &  67.7 &  40.3 &  62.9 &  28.6 &  57.6\\
 &    4 & 64.9 & 75.4 &   51 & 64.1 & 40.6 &  57.3 &  30.6 &  51.1 &  18.9 &  44.1\\
 &   10 & 54.8 & 64.2 & 37.3 & 50.1 & 23.4 &  41.1 &    14 &  32.8 &  5.85 &  24.1\\
 &  100 & 47.4 & 49.7 & 13.2 & 21.7 & 1.13 &  6.26 & 0.972 & 0.829 & 0.706 & 0.377\\
 & 1000 &   47 & 47.3 & 9.33 & 10.8 & 1.19 & 0.945 & 0.964 & 0.734 &  0.71 & 0.467\\
\hline \hline
\end{tabular}}
\caption{Empirical MISE  obtained with $1000$ samples and
approximations performed with $M=8$, for different sample  size
($n=100, 250, 500, 1000, 2500$) and different values of $s2n$ (2, 4,
10, 100, 1000, the higher $s2n$ the lower the noise level).}\label{BasicMise}
\end{table}
\end{center}

\begin{center}
\begin{table}
{\small
\begin{tabular}{||c|r||c|c||c|c||c|c||c|c||c|c||}\hline\hline
\multicolumn{2}{||c||}{}  & \multicolumn{2}{c||}{$n=100$}& \multicolumn{2}{c||}
{$n=250$}& \multicolumn{2}{c||}{$n=500$}& \multicolumn{2}{c||}{$n=1000$}&
\multicolumn{2}{c||}{$n=2500$}\\ \hline \hline
$\;\;\;g\;\;\;$ & $s2n$ & Lap. & Gaus. & Lap. & Gaus. & Lap. & Gaus. & Lap. & Gaus. & Lap. & Gaus. \\ \hline\hline
\multirow{3}{.2cm}{\rot{{\tiny Mix.Gam}}}
 &    2 & 0.96 & 0.92 & 1.3 & 1.1 & 1.4 & 1.6 & 1.3 & 2.2 & 1.3 & 1.6\\
 &    4 &    1 & 0.97 & 1.2 & 1.1 & 1.1 & 1.5 & 1.1 & 1.3 & 1.1 & 1.2\\
&   10 &    1 &    1 & 1.1 & 1.1 &   1 & 1.1 &   1 & 1.1 &   1 & 1.1\\
\hline\hline
\multirow{3}{.2cm}{\rot{{\tiny stable $\frac 12$}}}
 &    2 & 1.2 & 1.4 & 1.2 & 1.4 & 1.2 & 1.3 & 1.2 & 1.3 & 1.2 & 1.3\\
 &    4 & 1.1 & 1.2 & 1.1 & 1.2 & 1.1 & 1.2 & 1.1 & 1.1 & 1.1 & 1.1\\
  &   10 & 1.1 & 1.1 &   1 & 1.1 &   1 &   1 & 1.1 & 1.1 & 1.1 & 1.1\\
\hline\hline
\multirow{3}{.2cm}{\rot{{\tiny stable $\frac34$}}}
 &    2 & 1.8 & 2.2 & 1.6 & 2.2 & 1.6 & 2.1 & 1.6 &   2 & 1.6 & 1.8\\
 &    4 & 1.4 & 1.5 & 1.3 & 1.5 & 1.2 & 1.5 & 1.3 & 1.4 & 1.2 & 1.4\\
 &   10 & 1.2 & 1.2 & 1.1 & 1.2 & 1.1 & 1.1 & 1.1 & 1.2 & 1.1 & 1.2\\
\hline\hline
\multirow{3}{.2cm}{\rot{{\tiny Cauchy}}}
 &    2 & 1.5 & 2.4 & 1.4 & 2.3 & 1.2 & 2.1 &  1.1 & 1.9 & 0.86 &  1.5\\
 &    4 & 1.3 & 1.6 & 1.1 & 1.6 & 1.1 & 1.5 & 0.97 & 1.3 & 0.88 &  1.1\\
 &   10 & 1.1 & 1.2 &   1 & 1.1 &   1 & 1.1 & 0.96 &   1 & 0.91 & 0.96\\
\hline\hline
\multirow{3}{.2cm}{\rot{{\tiny Mix.Gau}}}
 &    2 & 1.1 & 2.1 &    1 & 1.7 &    1 & 1.5 & 0.97 & 1.2 & 0.94 &  1.1\\
 &    4 &   1 & 1.3 &    1 & 1.2 & 0.99 & 1.1 & 0.99 &   1 &    1 & 0.98\\
 &   10 &   1 & 1.1 & 0.99 &   1 & 0.99 &   1 & 0.99 &   1 &  1.1 &    1\\
\hline\hline
\multirow{3}{.2cm}{\rot{{\tiny F\'ejer 1}}}
 &    2 &  1.2 &  3.8 & 0.58 &  3.3 & 0.59 &  1.1 & 0.59 & 0.36 & 0.59 & 0.35\\
 &    4 & 0.87 &  2.3 & 0.81 &  1.1 & 0.81 & 0.73 &  0.8 & 0.73 & 0.83 & 0.71\\
 &   10 & 0.93 &  1.1 & 0.92 & 0.91 & 0.92 & 0.92 & 0.92 & 0.92 &  0.9 &  0.9\\
\hline \hline
\multirow{3}{.2cm}{\rot{{\tiny F\'ejer 5}}}
 &    2 & 1.4 & 1.5 & 1.6 & 1.7 &  1.8 &  1.8 &    2 &  1.8 &  2.2 &  1.8\\
 &    4 & 1.2 & 1.3 & 1.3 & 1.4 &  1.4 &  1.4 &  1.6 &  1.5 &  0.8 &  1.6\\
 &   10 & 1.1 & 1.1 & 1.2 & 1.2 &  1.3 &  1.3 & 0.92 &  1.4 & 0.94 & 0.87\\
\hline\hline
\multirow{3}{.2cm}{\rot{{\tiny F\'ejer 10}}}
 &    2 &  1.1 & 1.1& 1.1 & 1.2 &  1.2 & 1.2 &  1.3 &  1.2 &  1.5 &  1.2\\
 &    4 & 0.97 &    0.99 & 1.1 & 1.1 &  1.1 & 1.1 &  1.2 &  1.1 &  1.2 &  1.1\\
 &   10 & 0.95 &    0.96 &   1 &   1 &  1.1 & 1.1 &  1.1 &  1.1 &  1.1 &  1.1\\
\hline\hline
\multirow{3}{.2cm}{\rot{{\tiny F\'ejer 13}}}
 &    2 &    1 &  1.1 &  1.1 &    3 & 1.1 & 1.1 &  1.2 & 1.1 & 1.3 &  1.1\\
 &    4 & 0.94 & 0.95 &    1 &    1 &   1 &   1 &  1.1 & 1.1 & 1.1 &  1.1\\
 &   10 &  0.9 & 0.93 & 0.99 & 0.99 &   1 &   1 &  1.1 & 1.1 & 1.1 &  1.1\\
\hline\hline
\end{tabular}}
\caption{Ratio of the MISE with estimated $s2n$ over MISE with known $s2n$
when not strictly equal to one.}\label{Ratios2n}
\end{table}
\end{center}

\begin{center}
\begin{table}
\begin{tabular}{||c|r||c|c||c|c||c|c||c|c||c|c||}\hline\hline
\multicolumn{2}{||c||}{$\;$}  & \multicolumn{2}{c||}{$n=100$}& \multicolumn{2}{c||}{$n=250$}&
\multicolumn{2}{c||}{$n=500$}& \multicolumn{2}{c||}{$n=1000$}& \multicolumn{2}{c||}{$n=2500$}\\
\hline \hline
$g$ & $s2n$ & Lap. & Gaus. & Lap. & Gaus. & Lap. & Gaus. & Lap. & Gaus. & Lap. & Gaus. \\
\hline\hline
\multirow{3}{.2cm}{\rot{{\tiny $a$=0.5}}}
 &  2 & 1.4 & 1.2 & 1.3 & 1.1 & 1.3 & 1.2 & 1.2 & 1.4 & 1.1 & 1.3\\
 &  4 & 1.6 & 1.3 & 1.4 & 1.6 & 1.2 & 1.4 & 1.2 & 1.2 & 1.2 & 1.2\\
 & 10 & 1.6 & 1.6 & 1.4 & 1.5 & 1.4 & 1.4 & 1.3 & 1.2 & 1.1 & 1.1\\
 \hline\hline
\multirow{3}{.2cm}{\rot{{\tiny $a$=0.75}}}
 &  2 & 2 & 1.8 &   2 & 1.7 & 1.7 & 1.4 & 1.6 & 1.5 & 1.3 & 1.5\\
 &  4 & 3 &   2 & 2.4 & 2.6 &   2 & 2.1 & 1.6 & 1.7 & 1.4 & 1.5\\
 & 10 & 3 & 3.3 & 2.4 & 2.9 & 2.2 & 2.4 & 1.8 & 1.8 & 1.5 & 1.7\\
\hline\hline
\multirow{3}{.2cm}{\rot{{\tiny $a$=0.8}}}
 &  2 & 2.5 &   2 & 2.3 & 2.1 & 1.9 & 1.7 & 1.7 & 1.6 & 1.3 & 1.5\\
 &  4 & 3.5 & 2.4 & 2.7 & 3.1 & 2.3 & 2.5 & 1.8 & 2.1 & 1.5 & 1.7\\
 & 10 & 3.8 & 4.2 & 2.9 & 3.8 & 2.6 & 2.9 & 2.1 & 2.3 & 1.6 & 1.9\\
 \hline\hline
\multirow{3}{.2cm}{\rot{{\tiny $a$=0.9}}}
 &  2 & 4.6 & 3.5 & 3.8 & 3.6 & 3.3 & 2.9 & 2.4 & 2.6 & 1.8 & 2.1\\
 &  4 & 6.5 & 4.2 &   5 &   6 & 4.1 & 4.7 & 3.2 & 3.7 & 2.4 & 2.8\\
 & 10 & 7.5 &   8 & 5.8 & 6.9 & 4.9 & 5.4 & 3.7 & 4.3 & 2.6 & 3.4\\
\hline\hline
\multirow{3}{.2cm}{\rot{{\tiny $a$=0.95}}}
 &  2 & 8.2 & 6.1 & 7.5 & 24 & 5.9 &   5 & 4.2 & 4.4 & 2.8 & 3.2\\
 &  4 &  12 & 7.8 & 9.7 & 11 & 7.5 & 9.2 & 5.8 & 6.6 &   4 & 4.7\\
 & 10 &  15 &  16 &  11 & 14 & 9.5 &  11 & 7.2 & 7.8 & 4.9 & 6.4\\
\hline\hline
\end{tabular}
\caption{Ratio of the MISE obtained with $1000$ samples
and different values of $a$ over the MISE in the independent Gaussian
case \textbf{(k)} .}\label{DependMise1}
\end{table}
\end{center}



\begin{center}
\begin{table}
{\small
\begin{tabular}{||c|r||c|c||c|c||c|c||c|c||c|c||}\hline\hline
\multicolumn{2}{||c||}{}  & \multicolumn{2}{c||}{$n=100$}&
\multicolumn{2}{c||}{$n=250$}&
\multicolumn{2}{c||}{$n=500$}& \multicolumn{2}{c||}{$n=1000$}& \multicolumn{2}{c||}{$n=2500$}\\
\hline \hline
$g$ & $s2n$ & Lap. & Gaus. & Lap. & Gaus. & Lap. & Gaus. & Lap. & Gaus. & Lap. & Gaus. \\
\hline\hline
\multirow{3}{.2cm}{\rot{{\tiny $a$=0.5}}}
 &  2 & 1.3 & 1.3 & 1.3 & 1.3 & 1.3 & 1.3 & 1.2 & 1.2 & 1.1 & 1.1\\
 &  4 & 1.5 & 1.3 & 1.4 & 1.4 & 1.4 & 1.4 & 1.3 & 1.2 & 1.2 & 1.2\\
 & 10 & 1.6 & 1.5 & 1.5 & 1.5 & 1.4 & 1.4 & 1.3 & 1.3 & 1.2 & 1.2\\
\hline\hline
\multirow{3}{.2cm}{\rot{{\tiny $a$=0.75}}}
 &  2 & 2.2 & 1.9 & 2.2 & 1.9 &   2 & 1.9 & 1.8 & 1.7 & 1.5 & 1.5\\
 &  4 & 2.6 &   2 & 2.5 & 2.1 & 2.3 & 2.1 & 1.9 & 1.8 & 1.5 & 1.5\\
 & 10 & 2.9 & 2.6 & 2.7 & 2.5 & 2.5 & 2.3 &   2 &   2 & 1.6 & 1.5\\
\hline\hline
\multirow{3}{.2cm}{\rot{{\tiny $a$=0.8}}}
 &  2 & 2.6 & 2.2 & 2.6 & 2.2 & 2.5 & 2.1 & 2.2 & 1.9 & 1.7 & 1.6\\
 &  4 & 3.1 & 2.4 &   3 & 2.5 & 2.6 & 2.5 & 2.3 & 2.2 & 1.8 & 1.7\\
 & 10 & 3.5 &   3 & 3.5 &   3 & 2.9 & 2.7 & 2.5 & 2.4 & 1.9 & 1.8\\
\hline\hline
\multirow{3}{.2cm}{\rot{{\tiny $a$=0.9}}}
 &  2 & 4.4 & 3.5 & 4.7 & 3.9 & 4.5 & 3.6 & 3.7 & 3.2 & 2.7 & 2.4\\
 &  4 & 5.9 &   4 & 5.6 & 4.6 & 5.4 & 4.6 & 4.2 & 3.7 & 2.9 & 2.7\\
 & 10 & 6.8 & 5.6 & 6.7 & 5.8 & 5.6 & 5.4 & 4.5 & 4.4 & 3.3 & 2.8\\
\hline\hline
\multirow{3}{.2cm}{\rot{{\tiny $a$=0.95}}}
 &  2 & 8.4 & 5.9 & 9.1 & 6.9 & 8.6 & 7.1 & 7.3 & 5.9 & 4.7 & 4.3\\
 &  4 &  11 & 6.9 &  12 & 8.4 &  11 & 8.6 & 8.1 & 7.3 & 5.4 & 4.8\\
 & 10 &  12 &  11 &  14 &  12 &  12 &  12 & 9.1 & 9.1 & 6.3 & 5.3\\
\hline\hline
 \end{tabular}}
\caption{Ratio of  MISE obtained with $1000$ samples and different values of $a$ over the MISE
in the independent mixed Gaussian case \textbf{(l)}.}\label{DependMise2}
\end{table}
\end{center}

\begin{center}
\begin{table}
\begin{tabular}{||c|l||c|c||c|c||}\hline\hline
\multicolumn{2}{||c||}{$\times 10^{-2}$} & \multicolumn{2}{c||}{$n=100$} & \multicolumn{2}{c||}{$n=250$}\\\hline \hline
density $g$ & method & $\varepsilon$ Lap. & $\varepsilon$ Gaus. & $\varepsilon$ Lap. & $\varepsilon$ Gaus. \\ \hline\hline
\multirow{4}{3cm}{\hspace{0.7cm}(e) or $\# 2$\\ \hspace{1cm}$\chi^2(3)$ \\\hspace{0.7cm}($s2n$=4)}
 & DG, lower median& 1.5 & 1.8 & --- & --- \\
 & DG,  higher median  & 1.8 & 2.2 & --- & --- \\ \cline{2-6}
& Proj.: median & 1.8 & 2.1 & --- & --- \\
& Proj.: mean & 1.9 & 2.2 & --- & --- \\  \hline\hline
\multirow{4}{3cm}{\hspace{0.7cm}(f) or $\# 6$\\
\hspace{0.4cm}Mix.Gamma \\\hspace{0.7cm}($s2n$=10)}
&  DG, lower median & --- & --- & 0.21 & 0.23 \\
& DG,  higher median  & --- & --- & 0.24 & 0.26 \\ \cline{2-6}
& Proj.: median  & --- & --- & 0.17 & 0.27 \\
& Proj., mean & --- & --- & 0.20 & 0.28\\\hline\hline
\multirow{4}{3cm}{\hspace{0.7cm}(k) or $\# 1$\\ \hspace{0.7cm}
Gauss\\\hspace{0.7cm}($s2n$=4)}
& DG, lower median &  0.71 & 0.80 & 0.41 & 0.51 \\
& DG, higher median & 1.1 & 1.2 & 0.59 & 0.72 \\ \cline{2-6}
& Proj.: median & 0.45 & 0.76 & 0.31 & 0.22 \\
& Proj.: mean & 0.65 & 0.84 & 0.42 & 0.31\\ \hline \hline
\multirow{4}{3cm}{\hspace{0.7cm}(l) or $\# 3$\\
\hspace{0.6cm}Mix.Gauss\\\hspace{0.7cm}($s2n$=4)}
& DG, lower median &  1.8 & 2.7 & 1.1 & 2.0 \\
& DG,  higher median  & 3.1 & 3.4 & 2.3 & 2.8 \\ \cline{2-6}
& Proj.: median   & 0.48 & 0.62 & 0.23 & 0.26 \\
& Proj.: mean & 0.56 & 0.67 & 0.27 & 0.30\\\hline \hline
\end{tabular}
\caption{Lower and higher Median ISE obtained by Delaigle and
Gijbels (2004) with four different strategies
of bandwidth selection in kernel estimation compared with median and mean for our penalized
projection estimator.}\label{DelGibcomp}
\end{table}
\end{center}

\begin{center}
\begin{table}
\begin{tabular}{||c|c||c|c|c|c||}\hline
$\times 10^{-2}$ & method & $n=50$ & $n=100$ & $n=500$ & $n=1000$ \\\hline
$g$ Gaussian & D. Kernel & 1.18 & 0.63 & 0.13 & 0.08 \\ \cline{2-6}
& Gauss. Ker. & 1.72 & 1.27 & 0.28 & 0.16 \\ \cline{2-6}
& sinc Ker. & 2.16 & 1.14 & 0.26 & 0.10 \\ \cline{2-6}
             & Proj.  &   0.84   & 0.53 & 0.18 & 0.12 \\ \hline\hline
$g$ F\'ejer 5& D. Kernel & 2.29 & 0.79 & 0.22 & 0.13 \\\cline{2-6}
 & Gauss. Ker. & 3.07 & 1.84 & 0.55 & 0.22 \\ \cline{2-6}
 & sinc Ker. & 3.92 & 1.87 & 0.55 & 0.23 \\ \cline{2-6}
             & Proj. & 6.74     & 3.93 & 0.32 & 0.27 \\ \hline\hline
$g$ Gamma(2,3/2) & D. Kernel & 2.70 & 1.48 & 0.52  & 0.27 \\\cline{2-6}
 & Gauss. Ker. & 2.77 & 2.09 & 0.61 & 0.31 \\ \cline{2-6}
 & sinc Ker. & 6.17 & 4.03 & 1.66 & 0.37 \\ \cline{2-6}
             & Proj. &  3.13  & 2.19 & 0.96  & 0.65 \\ \hline
\end{tabular}
\caption{MISE for our projection estimator (Proj.) with
Laplace penalty using $s2n=10000$ and for direct density estimation by kernel of Dalelane(2004),
with Gaussian kernel (D. Kernel) or with $\sin(x)/x$ kernel (sinc).}\label{DelalaneComp}
\end{table}
\end{center}

\begin{center}
\begin{table}
{\small
\begin{tabular}{||c|r||c|c||c|c||c|c||c|c||c|c||}\hline\hline
\multicolumn{2}{||c||}{}  & \multicolumn{2}{c||}{$n=100$}&
\multicolumn{2}{c||}{$n=250$}&
\multicolumn{2}{c||}{$n=500$}& \multicolumn{2}{c||}{$n=1000$}& \multicolumn{2}{c||}{$n=2500$}\\
\hline \hline
$g$ & $s2n$ & Lap. & Gaus. & Lap. & Gaus. & Lap. & Gaus. & Lap. & Gaus. & Lap. & Gaus. \\
\hline\hline
\multirow{5}{.5cm}{\rot{Exp.}}
 &    2 & 2.7 & 2.3 & 3.4 & 2.7 & 3.9 & 3.1 & 4.6 & 3.4 & 5.6 & 3.8\\
 &    4 &   3 & 2.5 & 3.8 & 3.1 & 4.5 & 3.6 & 5.3 & 4.1 & 6.5 & 4.7\\
 &   10 & 3.4 &   3 & 4.5 & 3.8 & 5.4 & 4.6 & 6.7 & 5.3 & 8.6 & 6.4\\
 &  100 & 3.7 & 3.7 &   5 & 4.9 & 6.6 & 6.2 & 9.8 & 8.6 &  15 &  12\\
 & 1000 & 3.8 & 3.7 &   5 &   5 & 6.7 & 6.7 &  11 &  10 &  19 &  18\\
\hline\hline
\multirow{5}{.5cm}{\rot{Laplace}}
 &    2 & 1.4 & 1.3 & 1.3 & 1.3 & 1.3 & 1.3 & 1.3 & 1.3 & 1.3 & 1.3\\
 &    4 & 1.3 & 1.2 & 1.3 & 1.2 & 1.2 & 1.2 & 1.3 & 1.2 & 1.3 & 1.2\\
 &   10 & 1.3 & 1.2 & 1.2 & 1.2 & 1.3 & 1.2 & 1.3 & 1.2 & 1.4 & 1.2\\
 &  100 & 1.3 & 1.3 & 1.2 & 1.2 & 1.3 & 1.2 & 1.4 & 1.3 & 1.7 & 1.5\\
 & 1000 & 1.3 & 1.3 & 1.2 & 1.2 & 1.2 & 1.2 & 1.4 & 1.4 & 1.8 & 1.7\\
\hline\hline
\multirow{5}{.5cm}{\rot{Chi2(3)}}
 &    2 &  12 &  15 &  11 &  13 & 9.2 &  12 & 7.8 &  11 & 6.1 & 9.9\\
 &    4 &  12 &  15 & 9.6 &  12 & 8.1 &  11 & 6.6 & 9.2 & 5.1 & 7.7\\
 &   10 &  10 &  12 & 8.1 & 9.8 & 6.5 & 8.2 & 5.1 & 6.7 & 3.4 & 5.2\\
 &  100 & 9.6 & 9.9 & 7.4 & 7.7 & 5.6 & 5.9 & 3.9 & 4.3 & 1.9 & 2.3\\
 & 1000 & 9.6 & 9.6 & 7.4 & 7.3 & 5.6 & 5.6 & 3.9 & 3.9 & 1.7 & 1.8\\
\hline\hline
\multirow{5}{.5cm}{\rot{Cauchy}}
 &    2 & 4.6 & 6.1 & 4.2 & 5.3 & 3.8 & 4.9 & 3.3 & 4.6 & 2.7 & 4.3\\
 &    4 & 4.6 & 6.5 &   4 & 5.6 & 3.6 & 4.9 & 3.1 & 4.4 & 2.5 & 3.8\\
 &   10 &   4 & 5.8 & 3.5 & 4.5 & 3.1 & 3.9 & 2.5 & 3.4 & 2.2 & 2.7\\
 &  100 & 3.5 & 3.7 & 3.3 & 3.4 & 2.8 & 2.9 & 2.3 & 2.5 &   2 & 2.1\\
 & 1000 & 3.5 & 3.5 & 3.3 & 3.3 & 2.8 & 2.9 & 2.3 & 2.4 &   2 &   2\\
\hline\hline
\end{tabular}}
\caption{Ratio of the MISE obtained by method (E2) over the MISE obtained with
method (E1).} \label{E1E2comp}
\end{table}
\end{center}
\newpage

\begin{center}
\begin{table}
{\small
\begin{tabular}{||c|r||c|c||c|c||c|c||c|c||c|c||}\hline\hline
\multicolumn{2}{||c||}{$\;$} & \multicolumn{2}{c||}{$n=100$}&
 \multicolumn{2}{c||}{$n=250$} & \multicolumn{2}{c||}{$n=500$}&
\multicolumn{2}{c||}{$n=1000$}& \multicolumn{2}{c||}{$n=2500$}
\\ \hline \hline
\multicolumn{2}{||c||}{Noise} & Lap. & Gaus. &  Lap. & Gaus. &  Lap. & Gaus. &  Lap. & Gaus.&  Lap. & Gaus.\\
\hline
\multicolumn{2}{||c||}{Penalty} & Gaus. &  Lap. & Gaus. &  Lap. & Gaus. &  Lap. & Gaus.&  Lap. & Gaus. & Lap.\\
\hline\hline
$g$ & $s2n$ & &  &   & &  &  & & &  & \\
\hline\hline
\multirow{3}{.5cm}{\rot{{\tiny Laplace}}}
 &  2 & 0.93 &  1.1 & 0.92 &  1.2 &  1.1 & 1.3 &    1 & 1.5 & 1.6 &   2\\
 &  4 & 0.97 &    1 & 0.96 &    1 & 0.96 & 1.1 & 0.98 & 1.2 & 1.1 & 1.5\\
 & 10 & 0.99 & 0.99 & 0.99 & 0.99 &    1 &   1 & 0.99 &   1 &   1 & 1.2\\
\hline\hline
\multirow{3}{.5cm}{\rot{{\tiny Mix.Gam.}}}
&  2 & 0.98 & 1.1 & 0.93 & 1 & 0.91 & 1.1 &    1 & 1.2 & 1.2 & 1.5\\
 &  4 & 0.99 &   1 &    1 & 1 & 0.98 &   1 & 0.99 & 1.1 &   1 & 1.2\\
 & 10 &    1 &   1 & 0.98 & 1 & 0.98 & 1.1 & 0.98 &   1 &   1 &   1\\
\hline\hline
\multirow{3}{.5cm}{\rot{{\tiny Cauchy}}}
&  2 & 1.1 & 0.98 & 1 & 0.93 & 1.1 & 0.91 & 1.2 &    1 & 1.5 & 1.2\\
 &  4 &   1 & 0.99 & 1 &    1 &   1 & 0.98 & 1.1 & 0.99 & 1.2 &   1\\
 & 10 &   1 &    1 & 1 & 0.98 & 1.1 & 0.98 &   1 & 0.98 &   1 &   1\\
\hline\hline
\multirow{3}{.5cm}{\rot{{\tiny Gauss}}}
  &  2 & 0.95 &    1 & 0.93 &  1.2 & 0.88 & 1.2 &  1.2 & 1.1 &  1.5 & 1.1\\
 &  4 & 0.96 &    1 & 0.96 &    1 &    1 &   1 & 0.95 &   1 &  1.1 &   1\\
 & 10 &    1 & 0.97 &    1 & 0.96 &    1 &   1 & 0.97 & 1.1 & 0.99 &   1\\
\hline\hline
\multirow{3}{.5cm}{\rot{{\tiny F\'ejer 1}}}
 &  2 & 0.91 & 0.98 &  1.1 & 0.97 &  1.1 &    1 & 1.2 &  1.1 & 1.1 &  1.1\\
 &  4 & 0.99 &    1 &    1 &    1 & 0.95 & 0.99 &   1 & 0.96 &   1 &    1\\
 & 10 &    1 & 0.97 & 0.93 &    1 & 0.96 & 0.96 & 1.1 &    1 &   1 & 0.98\\
\hline \hline
 \end{tabular}}
\caption{Ratio between MISE with misspecified error
density (Laplace errors, $g$ estimated as if errors were Gaussian and reciprocally)
and MISE with correctly specified error density.}\label{errorinerror}
\end{table}
\end{center}

\begin{center}
\begin{table}
{\small
\begin{tabular}{||c|r||c|c||c|c||c|c||c|c||c|c||}\hline\hline
\multicolumn{2}{||c||}{$\;$}  & \multicolumn{2}{c||}{$n=100$}& \multicolumn{2}{c||}{$n=250$}&
\multicolumn{2}{c||}{$n=500$}&
\multicolumn{2}{c||}{$n=1000$}& \multicolumn{2}{c||}{$n=2500$}
\\ \hline \hline
$g$ & $s2n$ &  Lap. & Gaus.  &  Lap. & Gaus. &  Lap. & Gaus. &  Lap. & Gaus. &  Lap. & Gaus.   \\
\hline\hline
\multirow{3}{.5cm}{\rot{{\tiny Laplace}}}
 &  2 &    1 &  0.9 &  1.3 &  1.2 & 1.5 &  1.4 & 1.9 &  1.8 & 2.9 &
2.2\\
  &  4 & 0.95 & 0.68 &    1 & 0.87 & 1.2 &    1 & 1.5 &  1.3 & 2.3 &
1.9\\
  & 10 & 0.96 & 0.78 & 0.98 & 0.79 &   1 & 0.83 & 1.1 & 0.99 & 1.6 &
1.4\\
\hline\hline
\multirow{3}{.5cm}{\rot{\small{{\tiny Mix.Gam.}}}}
  &  2 &  0.9 & 0.89 & 0.92 & 0.78 & 0.99 & 0.81 & 1.2 &  1.1 & 1.8 &
1.6\\
  &  4 & 0.93 & 0.94 &  0.9 & 0.62 & 0.95 & 0.67 & 1.1 & 0.91 & 1.5 &
1.3\\
  & 10 & 0.96 & 0.89 & 0.92 & 0.65 & 0.96 & 0.77 &   1 & 0.92 & 1.2 &
1\\
\hline\hline
\multirow{3}{.5cm}{\rot{{\tiny Cauchy}}}
&  2 & 0.83 & 0.7 & 0.99 & 0.88 &  1.2 &  1.1 &  1.5 &  1.6 & 2.2 & 2.3\\
 &  4 & 0.81 & 0.5 & 0.89 & 0.71 & 0.99 & 0.93 &  1.2 &  1.2 & 1.7 & 1.9\\
 & 10 & 0.87 & 0.6 & 0.89 & 0.82 & 0.91 & 0.84 & 0.92 & 0.99 &   1 & 1.4\\
\hline\hline
\multirow{3}{.5cm}{\rot{{\tiny Gauss}}}
 &  2 &  1.2 &  1.2 &  1.6 &   2 &  1.8 & 2.4 &  2.1 & 3.1 & 2.8 & 4.4\\
  &  4 &  1.1 & 0.95 &  1.1 & 1.6 &  1.2 & 1.8 &  1.4 &   2 & 1.9 &   3\\
  & 10 & 0.94 &  1.1 & 0.87 & 1.1 & 0.87 & 1.1 & 0.87 & 1.1 &   1 & 1.6\\
\hline\hline
\multirow{3}{.5cm}{\rot{{\tiny F\'ejer 1}}}
  &  2 & 0.97 & 0.92 &  1.1 &  1.6 &  1.3 &  1.5 &  1.5 & 1.9 &    2 &
2.6\\
  &  4 & 0.96 & 0.82 & 0.97 &  1.4 & 0.99 &  1.3 &  1.1 & 1.5 &  1.4 &
2\\
  & 10 &  0.9 &  1.2 & 0.86 & 0.99 & 0.83 & 0.98 & 0.82 & 1.1 & 0.89 &
1.2\\
\hline\hline
\end{tabular}}
\caption{Ratio between MISE when ignoring noise and MISE with
correctly specified error density.}\label{noerrorwitherror}
\end{table}
\end{center}

\end{document}